\documentclass[10pt, conference, letterpaper]{IEEEtran}

\usepackage{graphicx}
\usepackage{blindtext}
\usepackage[nottoc]{tocbibind}
\usepackage{todonotes}
\usepackage{amsmath}
\usepackage{soul}
\usepackage{subcaption}
\usepackage{url}

\title{Science DMZ Networks: How Different are They Really?}

\author{
\IEEEauthorblockN{Emily Mutter and Susmit Shannigrahi}
\IEEEauthorblockA{Computer Science Department, Tennessee Tech University \\
Cookeville, TN, USA \\
\{ebmutter42, sshannigrahi\}@tntech.edu}
}

\begin{document}



\maketitle
\begin{tikzpicture}[remember picture, overlay]
    \node[anchor=north east, xshift=-2.5cm, yshift=-1cm, inner sep=5pt, fill=green, text=red, rotate=0] at (current page.north east) {This is the accepted version of the paper. The final version will appear in the proceedings of IEEE LCN 2024.};
\end{tikzpicture}
 \nopagebreak
\begin{abstract}


The Science Demilitarized Zone (Science DMZ) is a network environment optimized for scientific applications. A Science DMZ provides an environment mostly free from competing traffic flows and complex security middleware such as firewalls or intrusion detection systems that often impede data transfer performance. The Science DMZ model provides a reference set of network design patterns, tuned hosts and protocol stacks dedicated to large data transfers and streamlined security postures that significantly improve data transfer performance, accelerating scientific collaborations and discovery.

Over the past decade, many universities and organizations have adopted this model for their research computing. Despite becoming increasingly popular, there is a lack of quantitative studies comparing such a specialized network to conventional production networks regarding network characteristics and data transfer performance. We strive to answer the following research questions in this study: Does a Science DMZ exhibit significantly different behavior than a general-purpose campus network? Does it improve application performance compared to such general-purpose networks? Through a two-year-long quantitative network measurement study, we find that a Science DMZ exhibits lower latency, higher throughput, and lower jitter behaviors. However, we also see several non-intuitive results. For example, a DMZ may take a longer route to external destinations and experience higher latency than the campus network. While the DMZ model benefits researchers, the benefits are not automatic - careful network tuning based on specific use cases is required to realize the full potential of such infrastructure. 

\end{abstract}

\section{Introduction}

Science and engineering applications are generating data at an unprecedented rate, producing hundreds of Terabytes to Petabytes of data within a very short time. Additionally, scientific collaborations are becoming increasingly global, which means the researchers must transfer these datasets over the wide area networks to various scientific facilities. Such data transfers can occur between instruments, storage servers, computing systems, and cloud computing platforms. General-purpose enterprise networks are often unsuitable for these types of data transfers since these networks prioritize general usability and security over performance. Scientific data transfers can face several challenges, such as bandwidth throttling, packet loss, slow throughput due to firewalls, intrusion detection systems, and other middleboxes, resulting in lower throughput, higher latency, and increased jitter and packet loss\cite{kissel2013efficient}\cite{dart2013science}\cite{lee2021scalable}. These challenges ultimately result in lower scientific productivity. 

Organizations often tailor a portion of their network for scientific data transfers to address these challenges. Such a network is generally called a Science DMZ. Science DMZs prioritize data transfer performance through streamlined security postures, such as simple rule-based access control lists rather than stateful firewalls, and network tuning,  such as large Ethernet frames and larger TCP windows.



Science DMZ networks are widely deployed at US academic campuses and other countries. By the latest count, more than 200 Science DMZs\cite{nsfDMZ} are in the US alone. While they are widely deployed, there is a lack of comparative, quantitative studies on how Science DMZ networks differ from their general-purpose counterparts. To address this gap, we have observed a general-purpose production network alongside a Science DMZ at a university campus over the past two years. We have deployed multiple measurement instruments in both networks and external facilities. We have used a number of standard network measurement tools (iperf3, ping, traceroute) and developed our own comparison software to measure network parameters such as RTT, Throughput, Jitter, and Packet loss. Externally, we have looked into network traffic to and from large cloud platforms (Google Cloud) and the RIPE Atlas measurement platform. 

These long-running measurements allowed us to understand the nuances in performance differences on both networks. We confirm that a Science DMZ generally provides a better environment for data-intensive research. However, such benefits are not automatic, and these networks may be susceptible to higher latency, packet loss, and longer paths. Therefore, careful network planning and optimization based on the requirements of specific use cases (e.g., bulk data vs. real-time) must be a part of such infrastructure.

\section{Background}

\subsection{General-purpose Networks vs. Science DMZs}
Campus networks are typically designed to serve large numbers of users and devices, support various applications (e.g., email, web browsing, and video), and provide security and quality of service \cite{8494729}. Campus networks are also equipped with firewalls to maintain network security that often takes precedence over quality of service\cite{8494729}. Because most general-purpose data flows are small (KBs-MBs) and have a short duration, moderate bandwidth, latency, and loss rates are usually sufficient for these flows. Most traditional applications on a campus network can adapt to the network’s bandwidth and are not overly sensitive to packet loss or jitter.

On the contrary, scientific data is often at terabyte- and petabyte-scale \cite{8494729}. When packet loss occurs during such transfers, TCP reduces throughput to levels where it can take days to complete a single data transfer \cite{dart2013science}. Energy Sciences Network (ESNet) developed the Science DMZ \cite{ESNET} architecture to address these issues and transfer scientific data faster. A Science DMZ is a portion of a network designed for high-performance scientific applications. It is often separated from the campus network either physically or logically \cite{8494729}. Science DMZs also have a different security posture than enterprise networks. Instead of using multi-layer firewalls as in enterprise networks, Science DMZs use simple stateless Access Control Lists (ACLs) that allow line-rate packet processing\cite{8494729}\cite{ESNET}. These steps decrease packet loss and congestion and increase throughput \cite{ESNET}. 
The Science DMZs are also often limited to specific (and vetted) users and devices, eliminating many of the threats on general-purpose networks and allowing Science DMZs to be equipped with more lenient security policies\cite{8494729}.

\subsection{State of Science DMZ Deployment}


The Science DMZ model, since its conception by Dart et al. \cite{dart2013science}, has seen widespread adoption and evolution, addressing the growing data-intensive demands of scientific research. There are currently more than 200 \cite{nsfDMZ} deployments across various organizations. The model's effectiveness in handling large-scale data transfers has been recognized across various scientific disciplines.
Peisert et al. \cite{peisert2017medical} discuss the implementation of medical science DMZs, providing a secure yet high-performance network environment crucial for handling sensitive medical data.
Gonzalez et al. \cite{gonzalez2017bigdata} and Liu et al. \cite{liu2017widearea} have explored the challenges and solutions in monitoring and optimizing data transfers over international research network connections. These studies underscore the importance of efficient data transfer protocols, as also highlighted by Kissel et al.\cite{kissel2013efficient}, to support the high-bandwidth requirements of global scientific collaborations.

The evolution of Science DMZs encompasses advancements in data rate management using machine learning \cite{caicedo2024machine}, scalable designs considering the nature of research traffic \cite{lee2021scalable}, and explicit feedback mechanisms for congestion control \cite{vega2023explicit}. Gegan et al. \cite{gegan2020anomaly} and Mazloum et al. \cite{mazloum2023enhancing} have contributed to enhancing security and measurement capabilities within Science DMZs and general purpose networks, addressing the critical need for secure data environments in the wake of increasing cybersecurity threats.

\subsection{Studies on Science DMZ Performance}
A few studies have looked at Science DMZ and application performance. 
A study by Crichigno et al. \cite{crichigno2018comprehensive} provides a comprehensive guide to a Science DMZ and describes some performance measurements. This study covers protocols and equipment essential for a high-performance Science DMZ, including router and switch configurations. This tutorial highlights performance evaluations using ESnet and a lab testbed, focusing on the effectiveness of router and switch equipment for large-scale data transfers. Additionally, it examines TCP attributes, their impact on network performance, the significance of specific data transfer tools and security measures in Science DMZs, and how such software and equipment can create bottlenecks. 
Another study by Lee et al. \cite{lee2021scalable} examines scientific research traffic on the Korea Research Environment Open Network, proposing a scalable Science DMZ design and an iterative greedy algorithm. The design enables cost-effective sharing of Data-Transfer Nodes (DTNs), crucial but expensive components in a Science DMZ. This approach significantly reduces capital expenditures (CAPEX) by up to 79\% compared to traditional models where each user has a dedicated DTN. 
Vega et al. \cite{vega2023explicit} shows that a P4-based controller that enhances data transfer rates can significantly improve network performance compared to non-dedicated Science DMZ cyberinfrastructure. Using emulated hosts, the study shows that their model can improve the flow completion time of large scientific data flows by an average of 21.7\%.
Calyam et al. \cite{prasad6785344} present a case study demonstrating the architecture's effectiveness in enhancing remote scientific collaboration and simplifying network management for High-Throughput Computing services.
In \cite{evalperf}, researchers studied the effect of the Science DMZ on network performance. They created three networks: one with no DMZ and no firewall, one with no DMZ and a firewall, and a DMZ scenario. 
They show that the DMZ scenario returns the overall best results compared to the no DMZ, no firewall, and no DMZ, no firewall scenarios. 
There have been several other studies on Science DMZ performance and specific tunings\cite{caicedo2024machine, lee2021scalable, abhinit2022science, vega2023explicit, gegan2020anomaly, mazloum2023enhancing, crichigno2021application, crichigno2021data}. However, these studies focused on particular aspects of a DMZ, such as data transfer performance and network tuning. As a result, studies that demonstrate quantitative improvements of a DMZ over general-purpose networks still need to be done.

\section{Measurement Infrastructure Setup}
In this study, we focus on analyzing the performance of two distinct networks and comparing them: the Science DMZ and the campus commodity network on our university campus. This section summarizes the tools and infrastructure we used for our study. 

\subsection{Measurement Tools and Infrastructure}

\textbf{RIPE Atlas}: RIPE Atlas is a global network with numerous servers, measurement devices, and virtual machines for network measurement. The RIPE Atlas network is a collection of ``probes" that conduct measurements and provide a real-time understanding of the condition of the Internet. Probes can conduct ping, traceroute, SSL/TLS, DNS, NTP, and HTTP measurements to select targets \cite{RIPEatlas}. We utilize RIPE Atlas to perform ping and traceroute measurements to and from servers on our campus.

\textbf{perfSONAR}: perfSONAR (performance Service-Oriented Network monitoring ARchitecture) is an open-source network measurement toolkit\cite{tierney2009perfsonar}. It provides many tools within one package to test and measure network performance. These tools include latency, throughput, trace, and disk-to-disk measurements. perfSONAR identifies areas of poor performance, by both location within the network and by a window of time in which they occur, and flags these problem spots. For this study, we created dedicated PerfSonar nodes and utilized publicly available ones. 

\textbf{Google Cloud}: Google Cloud is a platform that is traditionally not used for network measurements. However, in our case, it is evident that several science use cases are utilizing the Google Cloud for their computations. As such, we quantified the network parameters to and from the cloud.
\textbf{Standard tools:}
In addition to these distributed measurement platforms, we utilized several standard tools, such as ping\cite{ping} and traceroute\cite{traceroute}. Traceroute provides the option to use both UDP and ICMP, and we utilized both. For performance measurements, we utilized iPerf3\cite{zurawski2013perfsonar} - a command-line tool that measures the throughput between two IP endpoints. It also returns bandwidth, throughput, packet loss, and jitter from the tests.
Finally, tcpdump, libpcap and Wireshark allow packet capture and analysis of traffic traces. We utilized all three for our analyses. 

\subsection{Measurement Servers}

\begin{figure}[!ht]
    \centering
    \includegraphics[width=0.45\textwidth]{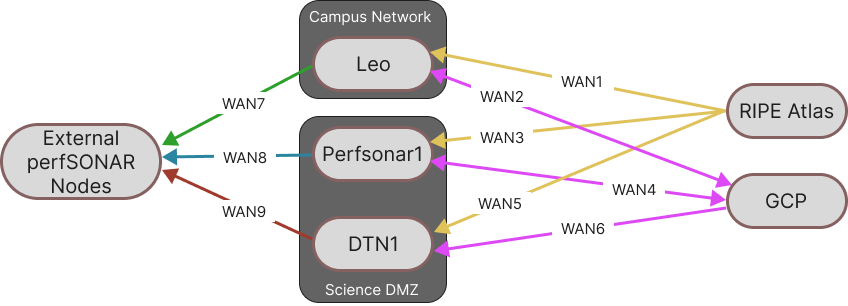}
    \caption{WAN Routes}
    \label{fig:wan}
\end{figure}

For these measurements, we created measurement servers within the campus network as well as on the DMZ. Figures \ref{fig:wan} and \ref{fig:lan} show these servers. The measurement server on the campus network is referred to as Leo. Leo ran a Perfsonar instance and had installed standard tools such as iperf3, ping, and traceroute.

On the Science DMZ, we used three other nodes: DTN1, DTN2, and Perfsonar1. We used DTN1 and DTN2 for data transfer experiments and Perfsonar1 for network measurement experiments. 

Externally, we used RIPE Atlas\cite{RIPEatlas}, publicly available PerfSonar nodes and Google Cloud (GCP) for our measurements.

\subsection{Network Routes}

A commercial ISP provided Layer3 network connectivity to the campus network. Internet2\cite{kratz2001ngi}, a network specifically designed to support scientific applications, provided Layer3 connectivity to the Science DMZ.

Internally, the campus network was connected to the provider using a 10Gbps link. All traffic passes through a gateway/firewall box that performs packet inspection. The Science DMZ network was connected to Internet2 at 10Gbps. This connection was served by a gateway and a security appliance using access control lists for security.

\begin{figure}[!ht]
    \centering
    \includegraphics[width=0.45\textwidth]{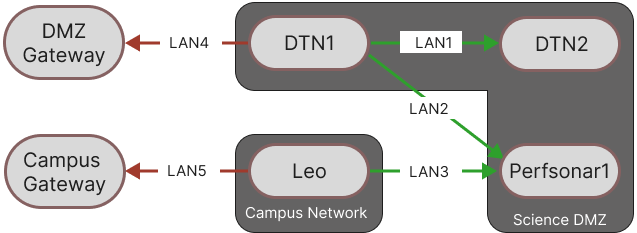}
    \caption{LAN Routes}
    \label{fig:lan}
\end{figure}

The campus and the Science DMZ network were logically separate. Even though they shared physical fibers, these networks used their own VLANs and traffic was completely separated. Figure \ref{fig:wan} shows the external routes. The colored lines in Figure \ref{fig:wan} show external (logical) connectivity to external measurement points (mainly RIPE Atlas and GCP). Figure \ref{fig:lan} shows local connections between Leo, DTNs, Perfsonar, and the gateways.


\subsection{Experiments}
 
\begin{table}[h!]
\centering
\caption{Measurement parameters for comparative analysis.}

\begin{tabular}{ |p{3cm}|p{4cm}|  }
 \hline
 \textbf{LAN-side Measurements} & \textbf{WAN-Side Measurements} \\ [0.5ex] 
 \hline
 \hline
 Throughput & Everything observed on LAN side  \\ 
 \hline
 RTT  between nodes & BGP routes to/from external vantage points \\
 \hline
 RTT between the node and the gateway & Path length between campus and external vantage points \\
 \hline
 Jitter & - \\
 \hline
 Packet loss & - \\
\hline
 \hline
\end{tabular}
\label{table:1}
\end{table}

We summarize our measurement experiments in Table \ref{table:1}. For this work, we conducted ``ping" tests to measure network latency, packet loss, and jitter. We utilized ``traceroute" to collect latency associated with network paths and identify intermediate hops between the source and destination nodes within each route. We utilized 
Iperf3 to observe throughput between external sources, the campus network, and the DMZ.

We originated these tests inbound from RIPE Atlas and Google Cloud Platform (GCP) virtual machines and outbound from the three local nodes (Leo, Perfsonar1, and DTN1).

\subsubsection{Internal clients $\rightarrow$ External servers Experiments}
Tests are run with one node of the campus network (Leo) and two nodes of the DMZ (Perfsonar1 and DTN1) posing as clients. 

We run ping and traceroute to 12 select perfSONAR nodes within the United States every 30 minutes. We send only ten packets during these tests so as not to overwhelm the external servers. We also used these clients to perform ping and traceroute from GCP VM instances hosted within the United States.

We used Iperf3 throughput experiments between two on-campus clients (Leo and Perfsonar1) and GCP VM instances, which were executed every 12 hours.

We perform the data transfer experiments using Leo and DTN1 as clients. We downloaded Linux ISOs from publicly available mirrors every four hours on both nodes and captured the packet headers using tcpdump.
These packet capture datasets allowed us to analyze interpacket delay, packet loss, round-trip time, packet retransmissions and download time. We observed the average daily value of these metrics in Wireshark, and we calculated the average RTT and interpacket delay externally and then plotted the daily values from DTN and Leo side-by-side.

\subsubsection{External Client $\rightarrow$ Internal Server Experiments}
We ran ping tests from RIPE Atlas to the local nodes every hour and traceroute tests every six hours. The ping and traceroute measurements send three packets of size 48 bytes during each execution. We executed a set of five experiments for each of these tests. For each experiment, we utilized five different RIPE Atlas source probes located within the United States.

Using the same method, we also run the ping and traceroute tests from GCP and the local nodes. Every 30 minutes, ping and traceroute tests run from GCP to the local nodes.

\subsubsection{Internal Clients $\leftrightarrow$ Internal Servers Experiments}
As previously mentioned, ping tests are performed to measure network latency, packet loss, and jitter, and traceroute tests are conducted to collect latency associated with network paths and identify intermediate hops between the source and destination nodes within each route. These tests are executed between the local network nodes (Leo, Perfsonar1, DTN1, DTN2), as well as between select local network nodes and the gateway to the campus network. 

Ping and traceroute tests are executed on these routes using the same method. Every 30 minutes, ping and traceroute tests run from Leo to DTN1, from Perfsonar1 to DTN1, from DTN1 to DTN2, and from both Leo and DTN1 to the gateway. Ping is designated to send only ten packets during the test.

\subsubsection{BGP Experiments}
For BGP experiments, we utilized a BGP dump from our Science DMZ BGP border router, which we manage. We obtained the BGP routes from our upstream provider on the campus network. 

\subsection{Data Analysis}
We parsed the collected data from ping, traceroute, and iperf3 into JSON and used Pandas, Seaborn, and Matplotlib to analyze and graph the results.

We examined the ping data to interpret latency, packet loss, and jitter. We analyzed the latency by taking all round-trip time (RTT) occurrences and graphing them with a Cumulative Distribution Function (CDF). We plotted daily packet loss by dividing the sum of all packets lost over a day by all packets sent over a day. We determined jitter by finding the difference in latency of subsequent packets. The jitter is then averaged daily and plotted with the standard deviation from that average.

We used traceroute data to calculate network latency and hop counts associated with network paths. We plot this by categorizing the measurements by the number of hops traversed in the network path and then averaging the latency observed for each route length. 

Finally, we used iperf3 and downloaded datasets for throughput insight. We plot this by averaging the bitrates from each day, categorizing them into ``sender" and "receiver," and then plotting the averages per day. 

\section{Results}

In this section, we discuss the comparative results from our experiments. We ran our experiments at regular intervals, as we described in the previous section. 

\subsection{Path Lengths}
\begin{figure*}[!h]
    \centering
    \def\imageA{images/EXTERN_LEO_TRACE.png}
    \def\imageB{images/EXTERN_PF1_TRACE.png}
    \def\imageC{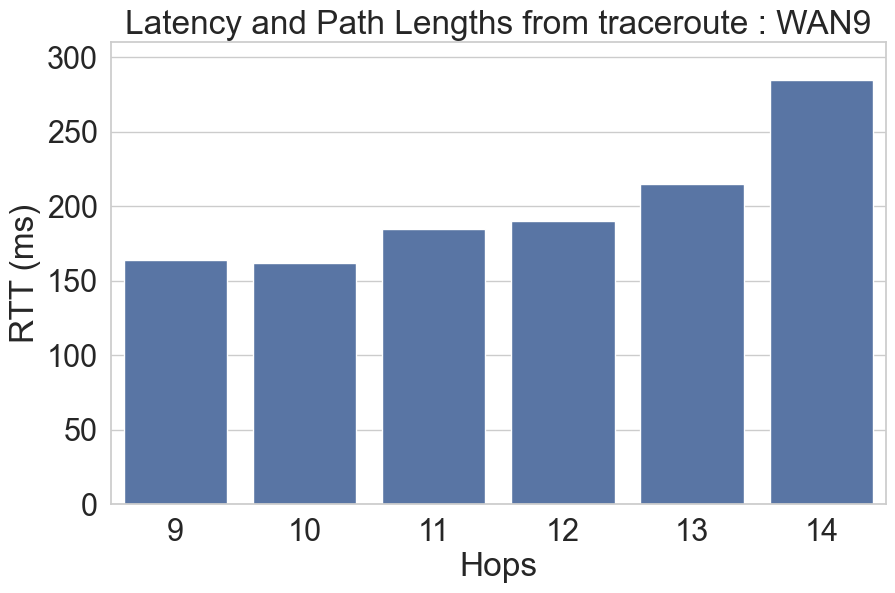}
    \def\imageD{images/GCP_LEO_TRACE.png} 
    \def\imageE{images/GCP_PF1_TRACE.png} 
    \def\imageF{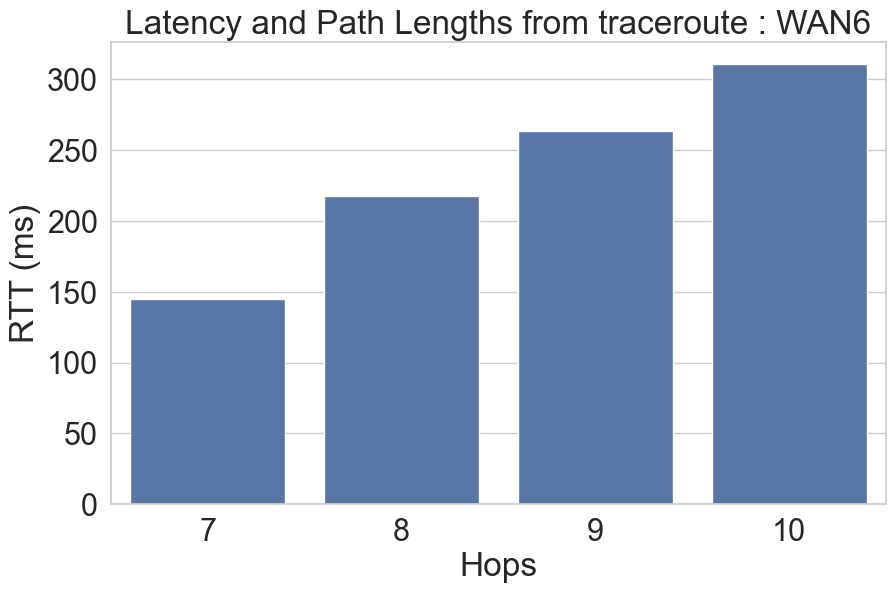} 
    \def\imageG{images/RIPE_LEO_TRACE.png} 
    \def\imageH{images/RIPE_PF1_TRACE.png} 
    \def\imageI{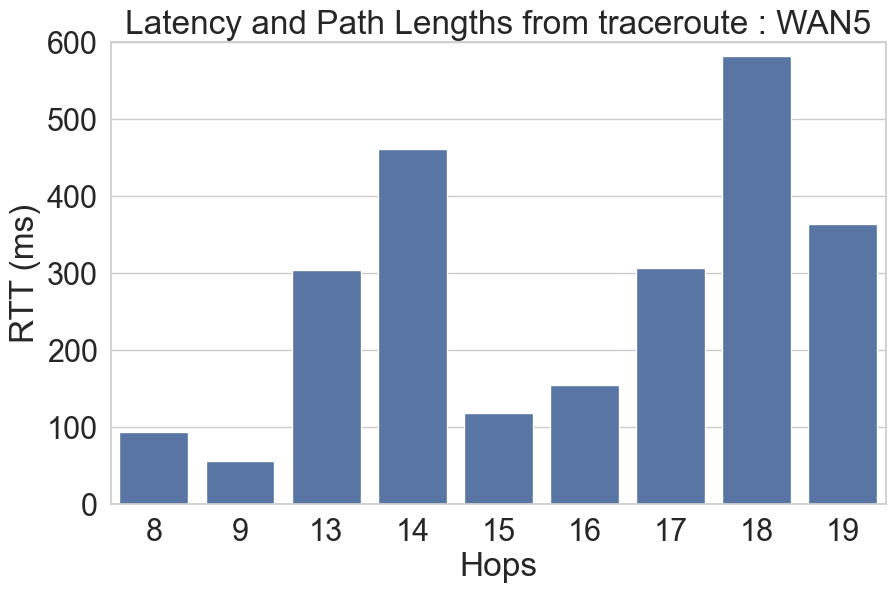} 
     \begin{subfigure}[t]{0.32\textwidth}
        \includegraphics[width=\columnwidth]{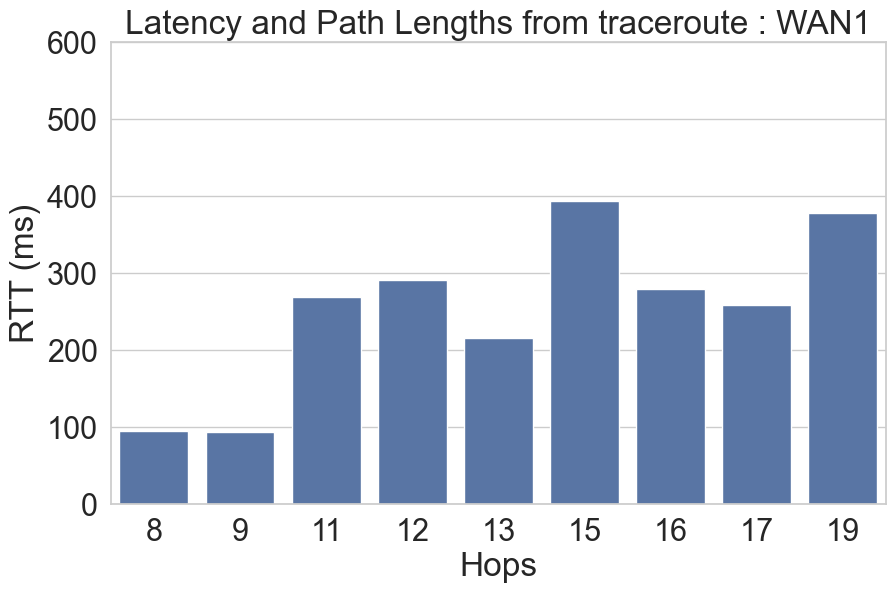}
        \caption{Latency and Path Lengths from RIPE Atlas to Leo}
        \label{fig:wan1_hops}
    \end{subfigure}
    \begin{subfigure}[t]{0.32\textwidth}
        \includegraphics[width=\columnwidth]{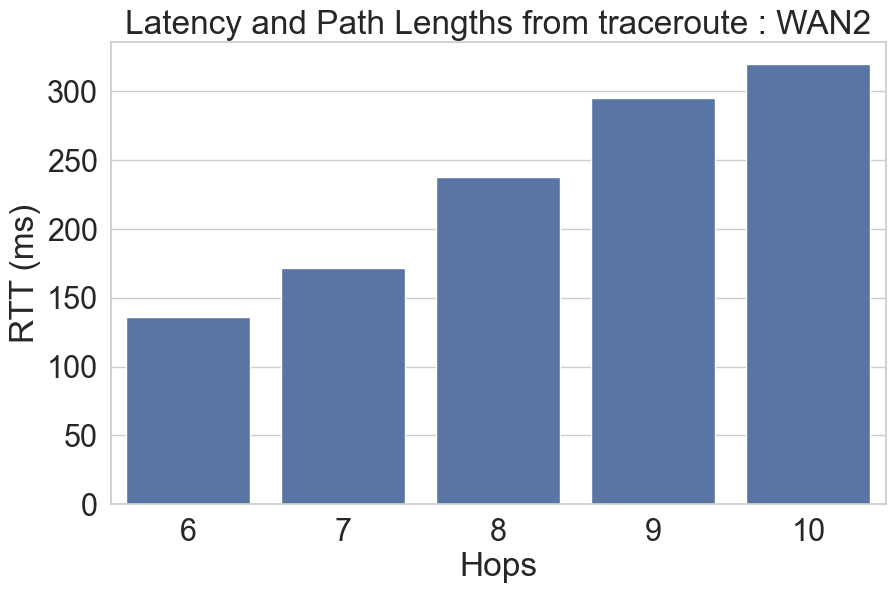}
        \caption{Latency and Path Lengths from GCP to Leo}
        \label{fig:wan2_hops}
    \end{subfigure}
    \begin{subfigure}[t]{0.32\textwidth}
        \includegraphics[width=\columnwidth]{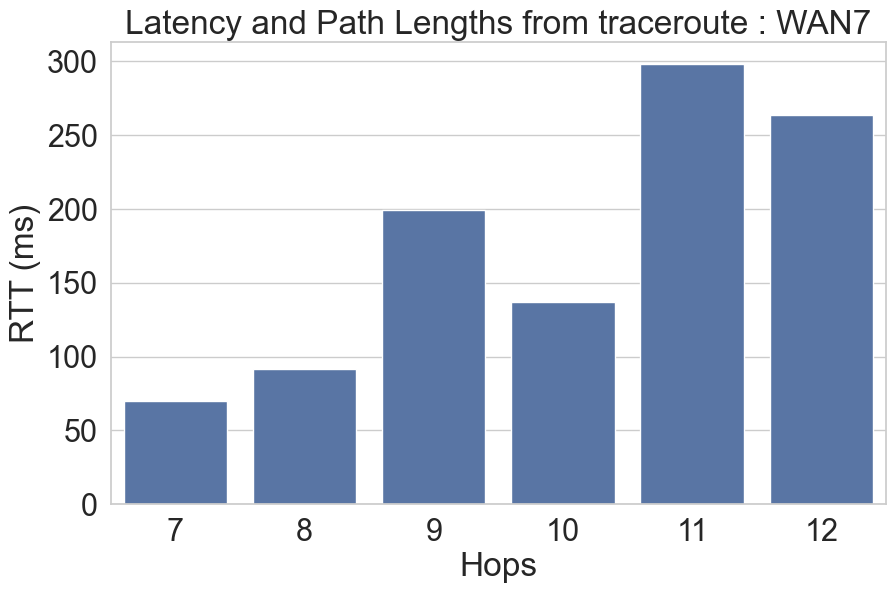}
        \caption{Latency and Path Lengths from Leo to external perfSONAR nodes}
        \label{fig:wan7_hops}
    \end{subfigure}
    \begin{subfigure}[t]{0.32\textwidth}
        \includegraphics[width=\columnwidth]{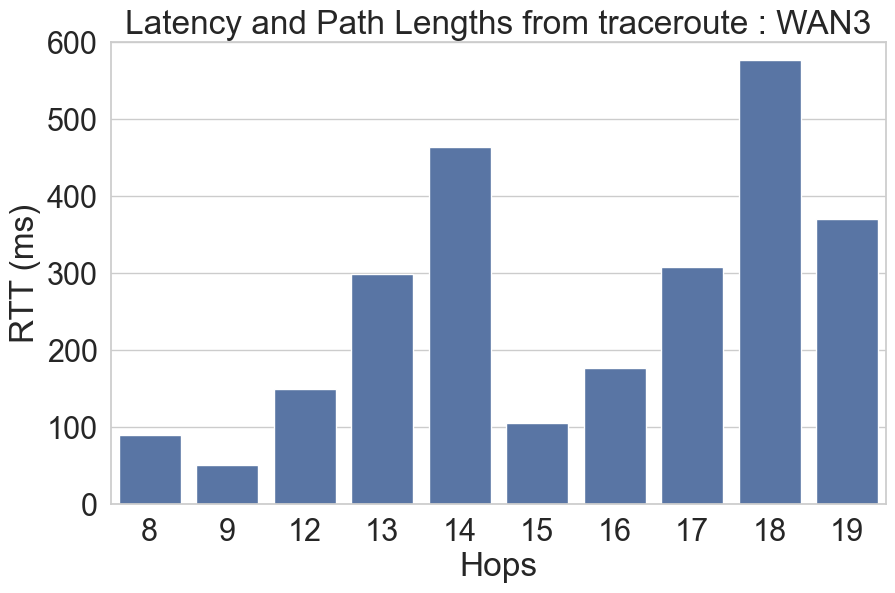}
        \caption{Latency and Path Lengths from RIPE Atlas to perfSONAR1}
        \label{fig:wan3_hops}
    \end{subfigure}
    \begin{subfigure}[t]{0.32\textwidth}
        \includegraphics[width=\columnwidth]{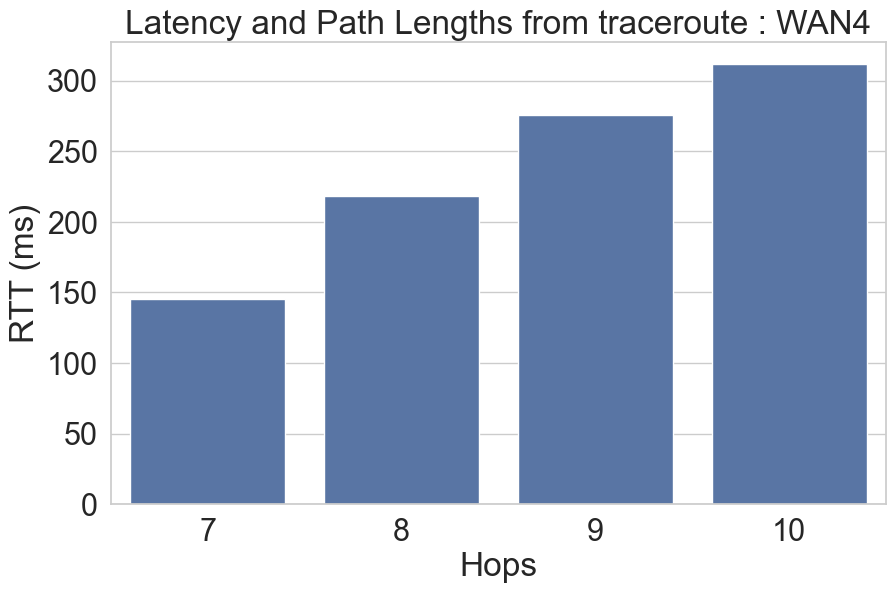}
        \caption{Latency and Path Lengths from GCP to perfSONAR1}
        \label{fig:wan4_hops}
    \end{subfigure}
    \begin{subfigure}[t]{0.32\textwidth}
        \includegraphics[width=\columnwidth]{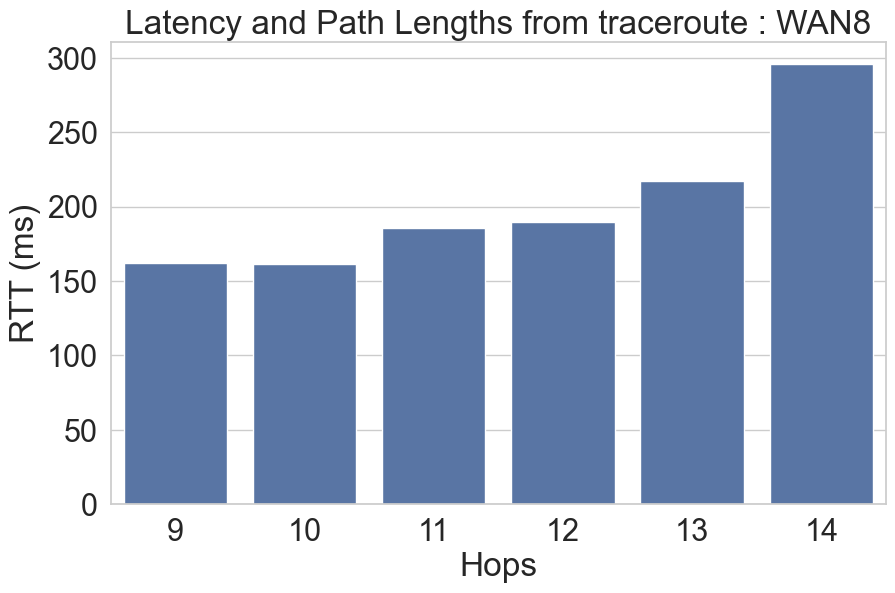}
        \caption{Latency and Path Lengths from perfSONAR1 to external perfSONAR nodes}
        \label{fig:wan8_hops}
    \end{subfigure}
\caption{Comparison of Path Lengths}
\label{fig:path_length_analysis}
\end{figure*}
Different upstream providers serve the DMZ and the commodity network in this study. A commercial ISP serves the campus network while the DMZ is served by Internet2, which is a specialized network for research. These experiments compare the path lengths of network destinations to/from internal and external vantage points.
Figures \ref{fig:wan1_hops} and \ref{fig:wan3_hops} show the average latency and path lengths between RIPE Atlas, Leo (located in the campus network), PerfSONAR1, and DTN1 (both located in the DMZ). In both experiments, The maximum hop counts are 19 hops, and the minimum is 8 hops. 

The latency and hop counts are lower between these servers and GCP, shown in Figures \ref{fig:wan2_hops} and \ref{fig:wan4_hops}. The hop count to these servers is 10 hops compared to 19 from RIPE Atlas. RIPE probes are hosted by various organizations and served by various ISPs. However, Google has a more optimized peering presence, leading to lower hop counts. The latency between GCP and these servers is also lower. Both for the DMZ and the campus network, the maximum latency is 300ms.
 But the DMZ exhibits lower latency at all route lengths in common with the campus network by $\sim$3\% - 6.78\%.
   
As exhibited in Figures \ref{fig:wan7_hops} and \ref{fig:wan8_hops}, when traffic is outbound to external perfSONAR nodes, Leo experiences routes with ranges 1-2 hops shorter than DMZ routes, and there is a point when the commodity network performs faster than the DMZ by 12.5\% at 10 hops. However, the DMZ tends to have a latency 20\% - 36.7\% lower than Leo, exclusively comparing common path lengths. Plots of the two DMZ nodes are very similar for this experiment, so Figure \ref{fig:wan8_hops} was selected to represent both nodes. However, we noticed one difference. The DTN1 node on the DMZ has a latency, at the longest path length of 14 hops, that is $\sim$6.75\% lower than that of the Perfsonar1 node on the DMZ. In these outbound experiments, the path lengths are between 7-12 hops on the campus network and 9-14 hops on the DMZ side. Since IP routing can be asymmetric, there is a mismatch between the hop counts from the inbound and the outbound experiments.

\textbf{Takeaways:} Given that a specialized research network serves the DMZ, Internet2, we expected this to have lower hop counts for inbound and outbound traffic. However, the DMZ experiments consistently show higher hop counts than the campus network. This observation is critical for delay-sensitive research applications, such as AR-VR, since moving them into the DMZ will potentially increase their hop count, resulting in end-to-end delay. 

From these experiments, we conclude that just placing research use cases into a DMZ may not automatically improve their performance/latency. Careful discussions and planning with upstream providers are needed to optimize routing and/or physical path. On our campus, we discovered the upstream provider routing traffic using a longer but less congested physical path rather than a short but more heavily used physical path.

\subsection{Latency}

Comparing the distribution in latency in Figures {3a-3f}, we find that the latencies are between 100-400ms on the campus network and 50-600ms on the DMZ. There is a significant spike in latency at the penultimate hop (DMZ gateway) for the DMZ experiments. 

More interestingly, the latencies are slightly higher on the DMZ for outbound experiments since the paths are typically longer. On the paths with higher hop counts, both the campus network and the DMZ experience similar latency as Figures \ref{fig:wan7_hops} and \ref{fig:wan8_hops} show.

Figures {4a-4c} compare the latency for inbound WAN traffic from RIPE Atlas and GCP and for outbound WAN traffic to perfSONAR nodes. The 95 percentile latency from RIPE Atlas to both the DMZ and campus network is around 80 ms. The 95 percentile latency from GCP to the campus network is around 35 ms, and to the DMZ is near 37 ms. This can again be attributed to better peering provided by GCP. 

Based on Figures \ref{fig:wan1_hops} and \ref{fig:wan3_hops}, when traffic is inbound from RIPE Atlas, the range of hops is the same to reach Leo and the two DMZ nodes. Comparing only standard path lengths, both nodes on the DMZ have similar latency, represented only by Figure \ref{fig:wan3_hops} for conciseness. A difference in latency was noted when the route averaged 16 hops for the DMZ nodes. At that point, the DTN1 node had 11\% lower latency. Both DMZ nodes often exhibit 34-73\% lower latency than Leo, but Leo has path lengths that have 13-30\% lower latency than the DMZ nodes.

\begin{figure*}[!htbp]
    \centering
    \def\imageA{images/EXTERN_PING_LATENCY.png}
    \def\imageB{images/GCP_PING_LATENCY.png}
    \def\imageC{images/RIPE_PING_LATENCY.png}
    \def\imageD{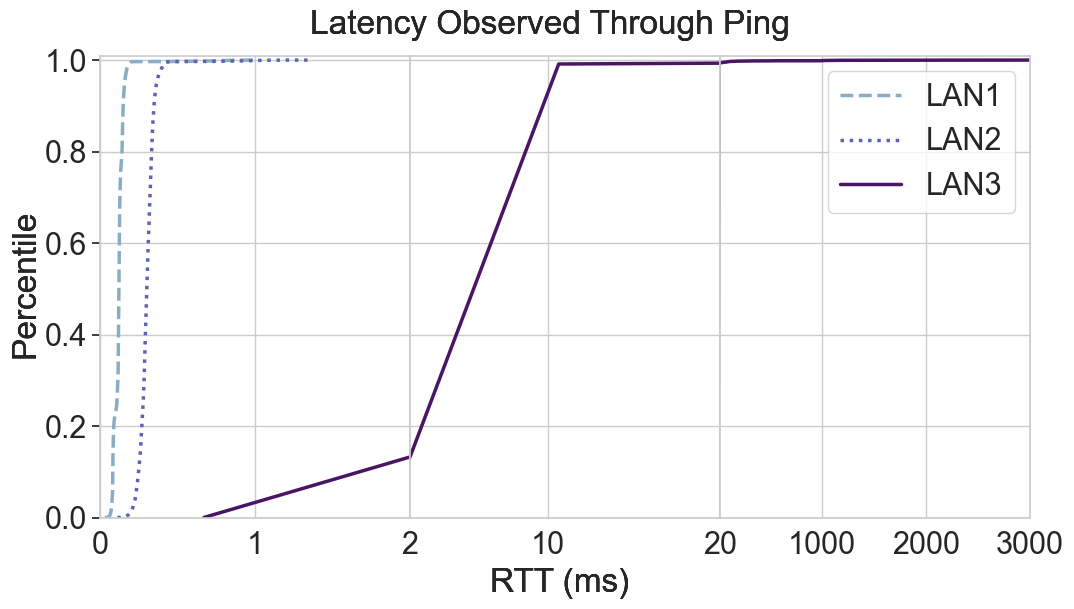} 
    \def\imageE{images/GATEWAY_PING_LATENCY.png} 

    \begin{subfigure}[t]{0.24\textwidth}
        \includegraphics[width=\linewidth]{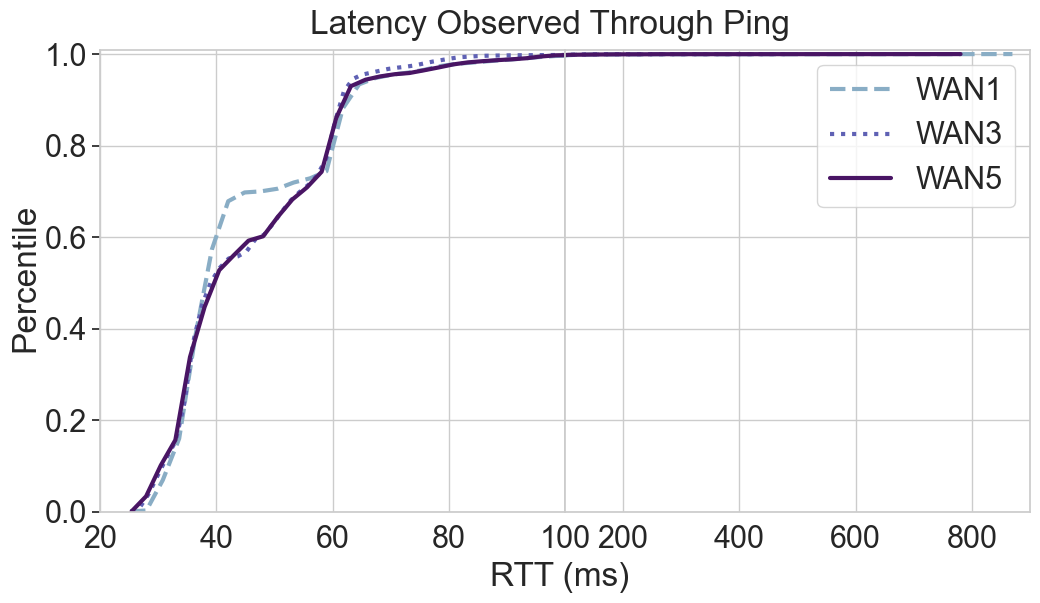}
        \caption{Latency from RIPE Atlas to campus and DMZ nodes}
        \label{fig:RIPE_latency}
    \end{subfigure}
    \hfill
    \begin{subfigure}[t]{0.24\textwidth}
        \includegraphics[width=\linewidth]{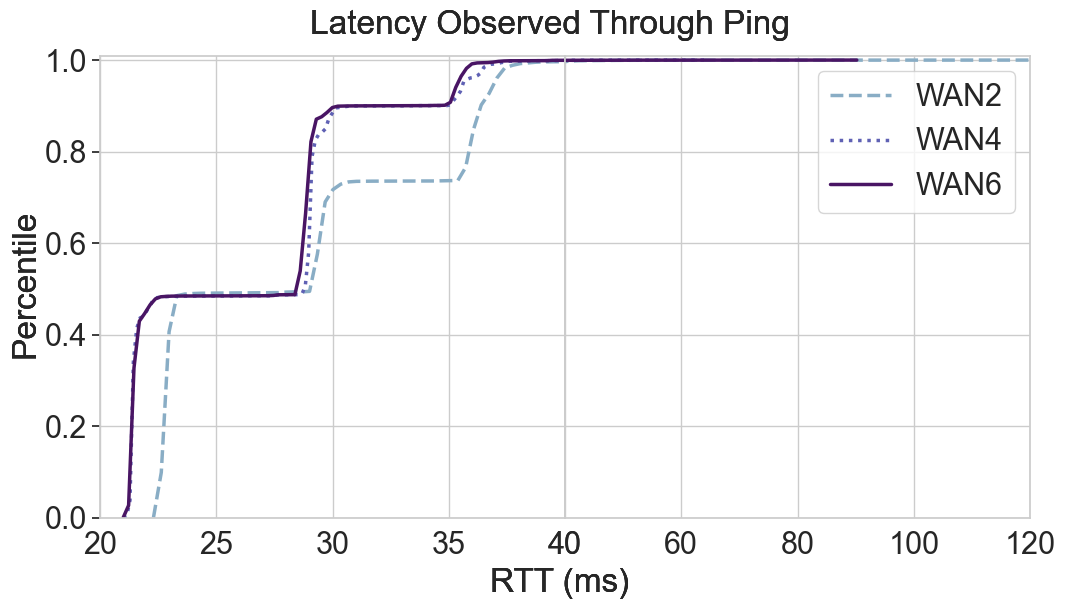}
        \caption{Latency from GCP to campus and DMZ nodes}
        \label{fig:GCP_latency}
    \end{subfigure}
    \hfill
   \begin{subfigure}[t]{0.24\textwidth}
        \includegraphics[width=\linewidth]{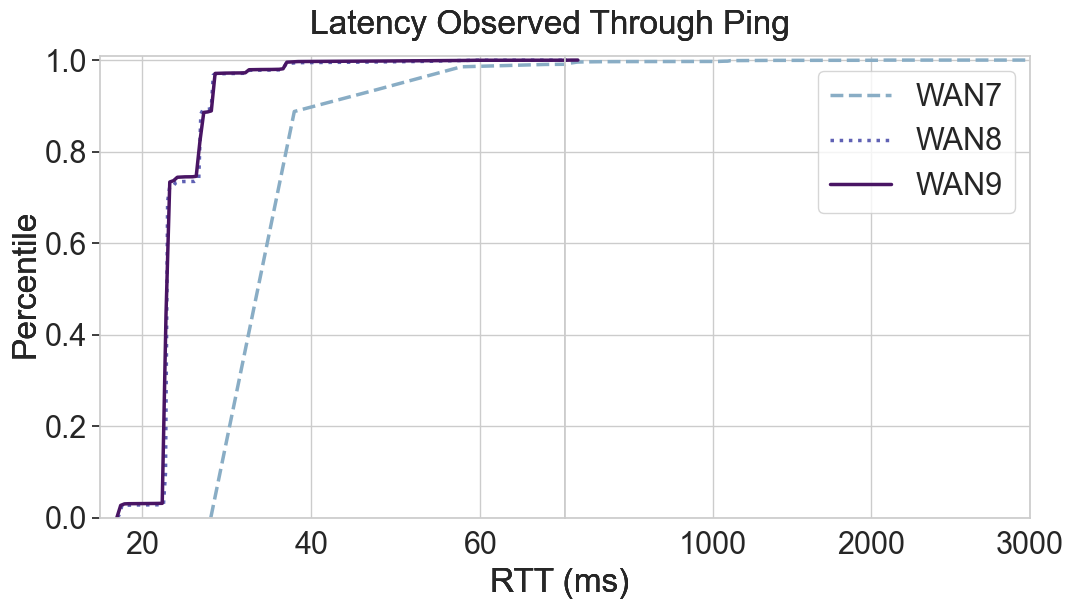}
        \caption{Latency from campus and DMZ to external perfSONAR nodes}
        \label{fig:perfSONAR_latency}
    \end{subfigure}
    \begin{subfigure}[t]{0.24\textwidth}
        \includegraphics[width=\linewidth]{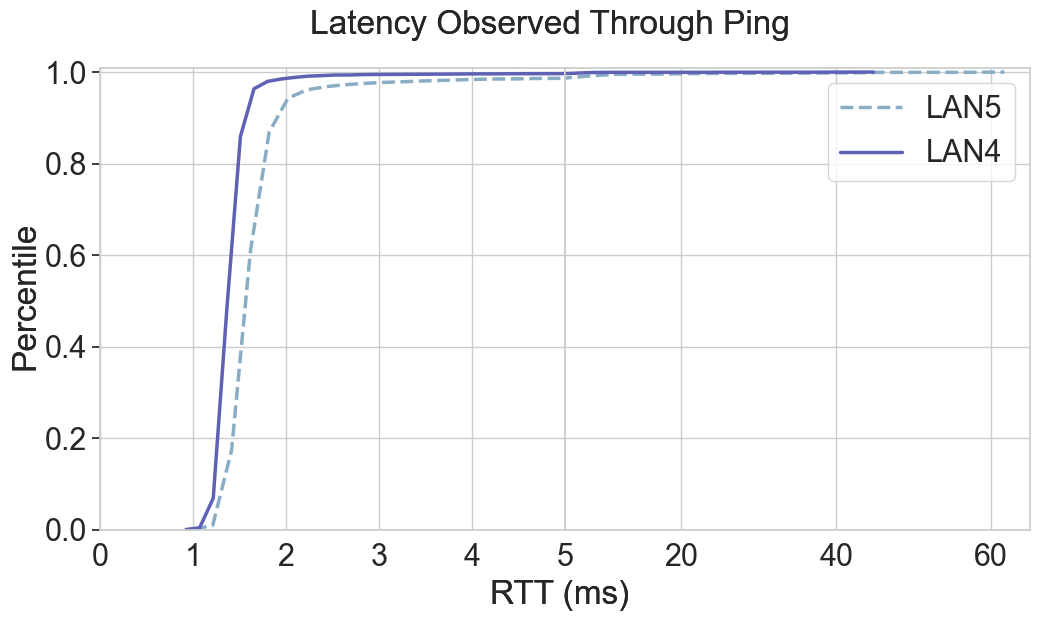}
        \caption{Latency from campus and DMZ nodes to campus gateway}
        \label{fig:gateway_latency}
    \end{subfigure}
    \caption{Latency Comparison}
    \label{fig:latency_analysis}
\end{figure*}

Based on Figures \ref{fig:wan2_hops} and \ref{fig:wan4_hops}, when traffic is inbound from Google Cloud, both nodes on the DMZ tend to have similar latency, with an occasional $\sim$2\% difference. Due to close similarities in their plots, only Figure \ref{fig:wan4_hops} represents the DMZ nodes for this experiment. The campus network tends to have similar latency to the DMZ or higher latency by $\sim$2\% - 24\%.

For the outbound experiments presented in Figure \ref{fig:perfSONAR_latency}, 95 percentile latency to external PerfSonar nodes is also around 35 ms on the DMZ side. On the campus network , the 95 percentile latency is near 55 ms. When traffic is outbound to perfSONAR nodes, both nodes on the DMZ exhibit similar latency, while the campus network experiences latency that is 30.43\% - 83\% slower.

Internally, we find the latency between the campus and DMZ nodes to be very low. However, given that the path length is minimal, the effect of the firewall is really pronounced here. Most pings between campus network servers and the DMZ exhibit a 10ms delay. Compare this to Figure \ref{fig:gateway_latency} where pings between the nodes and their respective gateways are less than 2ms. The inline firewall and access control lists (ACLs) add 8ms latency to each packet, which is very large. Most of these additional delays can be attributed to the firewall and packet inspection middleware.

\textbf{Takeaway:} We find that both the campus network and the DMZ exhibit similar latency, with the campus network occasionally having a lower average latency than the DMZ by as much as $\sim$20ms (5\% - 30.5\%). We find the measurements often get delayed on the DMZ (e.g., pings not arriving or arriving with higher latency), which affects results poorly. 
For internal measurements, we find that firewalls negatively affect performance, even when measurement boxes are placed on the same campus/data center.

\subsection{Packet Loss}
\begin{figure*}[!h]
    \centering
    \def\imageA{images/EXTERN_PACKET_LOSS.png}
    \def\imageB{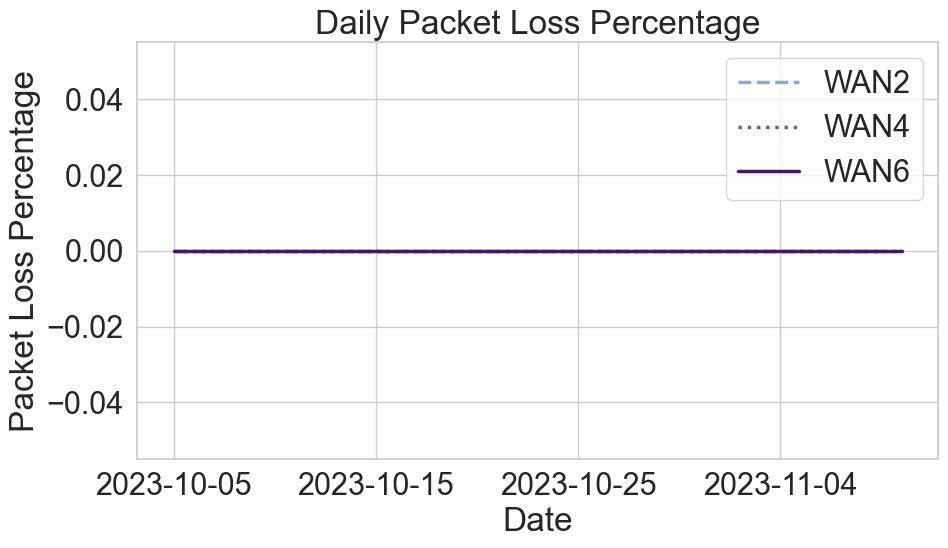}
    \def\imageC{images/RIPE_PACKET_LOSS.png}
    \def\imageD{images/LOCAL_PACKET_LOSS.png} 
    \begin{subfigure}[t]{0.324\textwidth}
        \includegraphics[width=\linewidth]{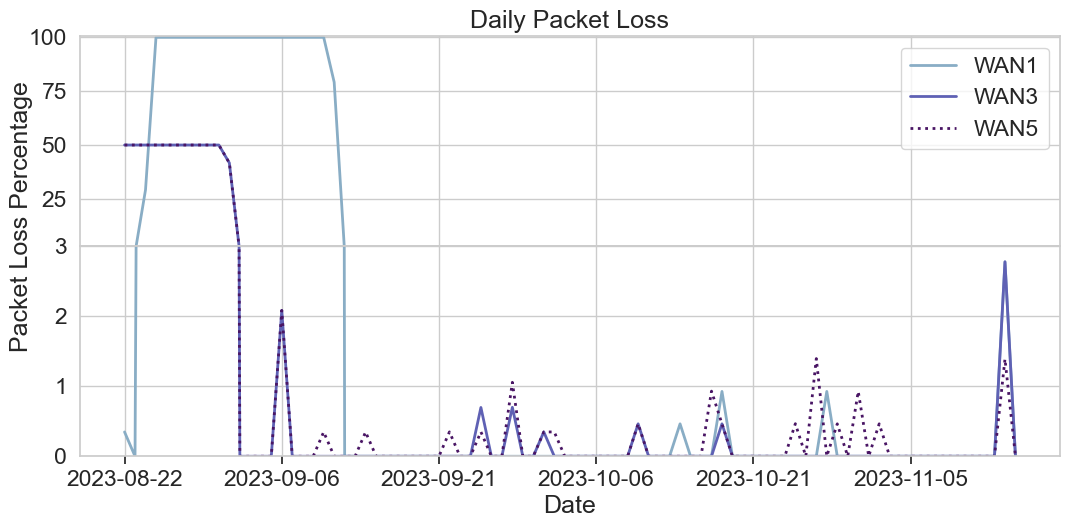}
        \caption{Packet loss from RIPE Atlas to campus and DMZ nodes}
        \label{fig:RIPE_loss}
    \end{subfigure}
   \begin{subfigure}[t]{0.324\textwidth}
        \includegraphics[width=\linewidth]{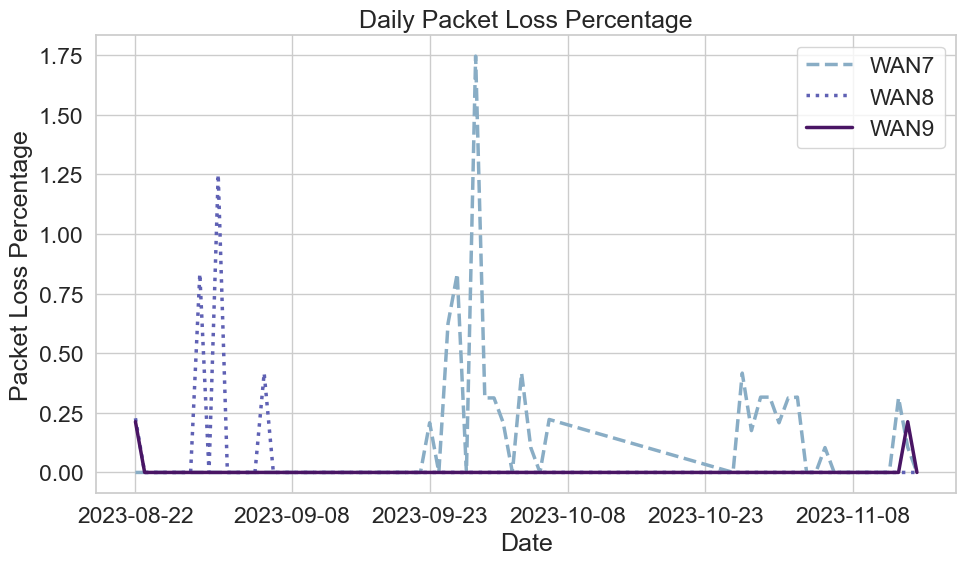}
        \caption{Packet loss from campus and DMZ to external perfSONAR nodes}
        \label{fig:perfSONAR_loss}
    \end{subfigure}
    \begin{subfigure}[t]{0.324\textwidth}
        \includegraphics[width=\linewidth]{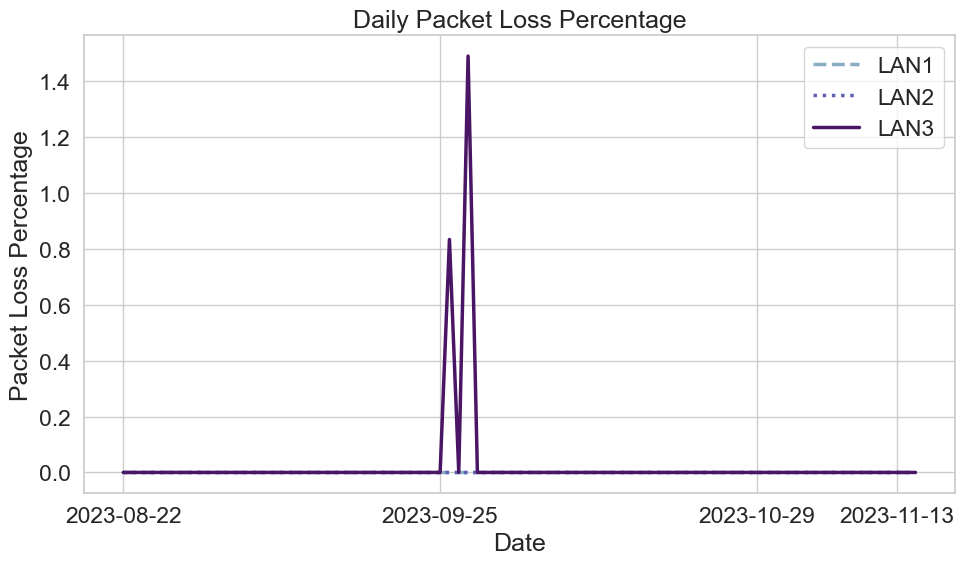}
        \caption{Packet loss between campus and DMZ nodes}
        \label{fig:internal_loss}
    \end{subfigure}
    \caption{Packet Loss comparison}
    \label{fig:packet_loss}

\end{figure*}

The DMZ experiences more packet loss than the campus network for inbound traffic from RIPE Atlas. While Leo exhibits a period of 100\% packet loss due to the campus node being down, as Figure \ref{fig:RIPE_loss} shows, both nodes of the DMZ experience 50\% genuine packet loss even when the network was up.
However, the packet loss is more consistent on the campus network, where we can observe regular 1-2\% packet losses. 

The Perfsonar1 node on the DMZ exhibits more packet loss than the DTN1 node on the DMZ; it loses $\sim$5\% more packets than DTN1 over three months, as Figure \ref{fig:RIPE_loss} shows. This is potentially because more experiments were conducted on the PerfSonar1 node than on DTN1.

When traffic is incoming from Google Cloud, there is no packet loss pattern across DMZ or campus network.
    
When traffic is outbound to external perfSONAR nodes, the campus network experiences more packet loss than the DMZ, but the Perfsonar1 node experiences more packet loss than the DTN1 node. Perfsonar1 exhibits $\sim$2\% more packet loss than DTN1. The campus network exhibits $\sim$.2\% more packet loss than Perfsonar1 and $\sim$5.7\% more packet loss than DTN1 as Figures \ref{fig:perfSONAR_loss} shows.

Internally, shown in Figure \ref{fig:internal_loss}, there is no pattern of packet loss sourced from the DMZ nodes, but there are spikes of loss sourced from the campus node of less than 2\%.

\textbf{Takeaways: } The campus network experiences more regular packet loss. Firewalls and middleboxes contribute to these packet loss events. Packet loss also occurs on the outbound paths from campus, again, potentially due to the presence of firewalls. This observation is important since large data transfers are sensitive to packet loss. Placing research use cases on a shared campus network will affect data transfer performance. Such use cases should be placed in a DMZ network, which has a lower loss rate due to the simplified nature of such networks.
\subsection{Jitter}
\begin{figure*}[!htbp]
    \centering
    \def\imageA{images/EXTERN_LEO_JITTER.png}
    \def\imageB{images/EXTERN_PF1_JITTER.png}
    \def\imageC{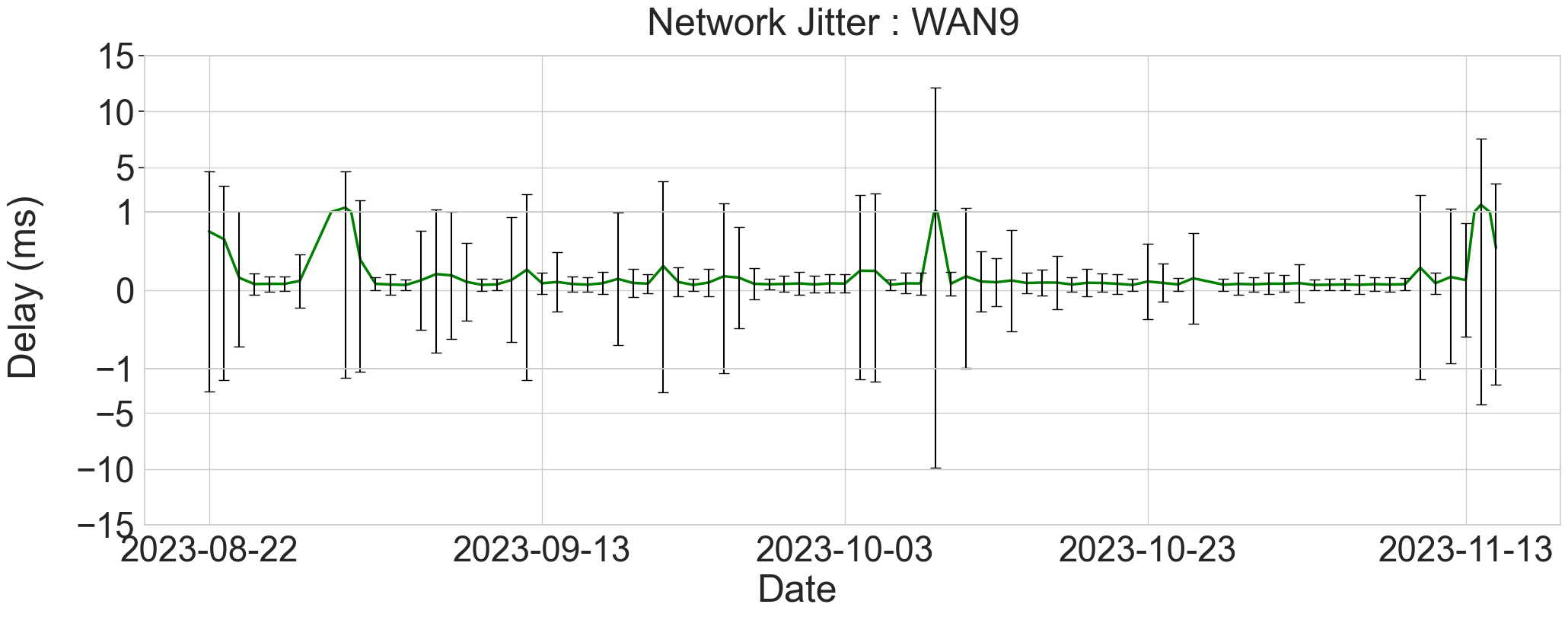}
    \def\imageD{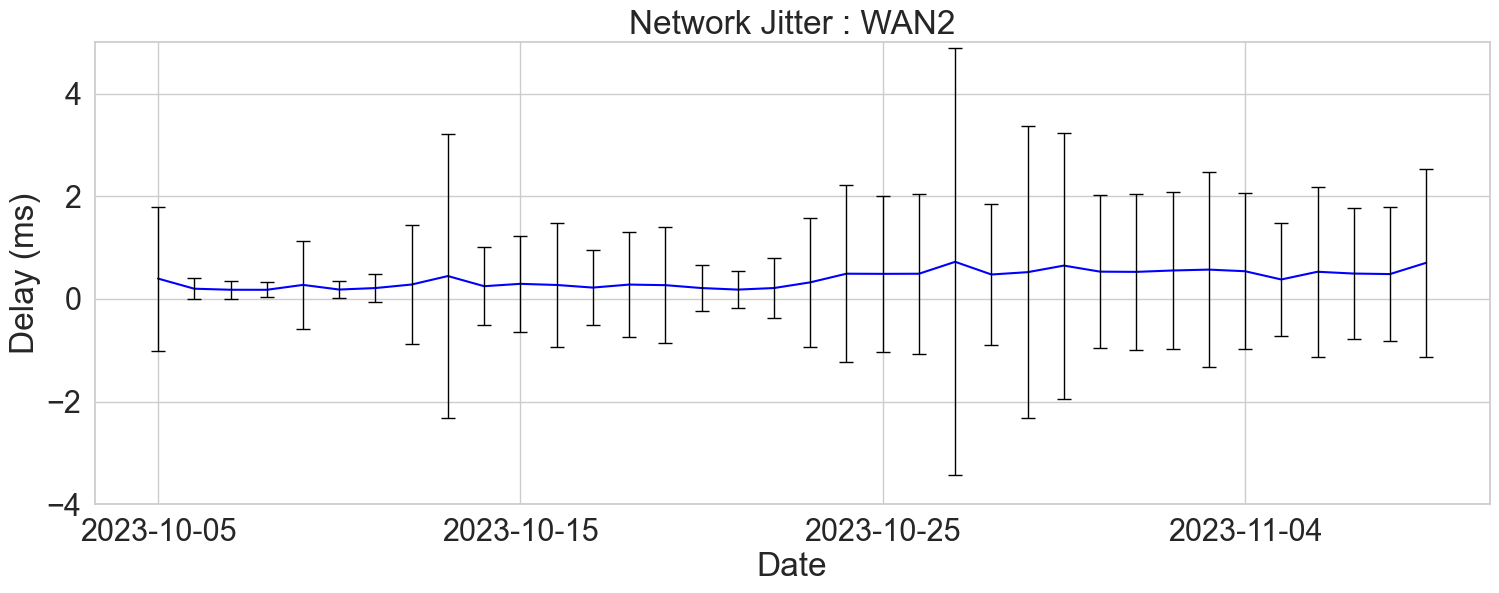} 
    \def\imageE{images/GCP_PF1_JITTER.png} 
    \def\imageF{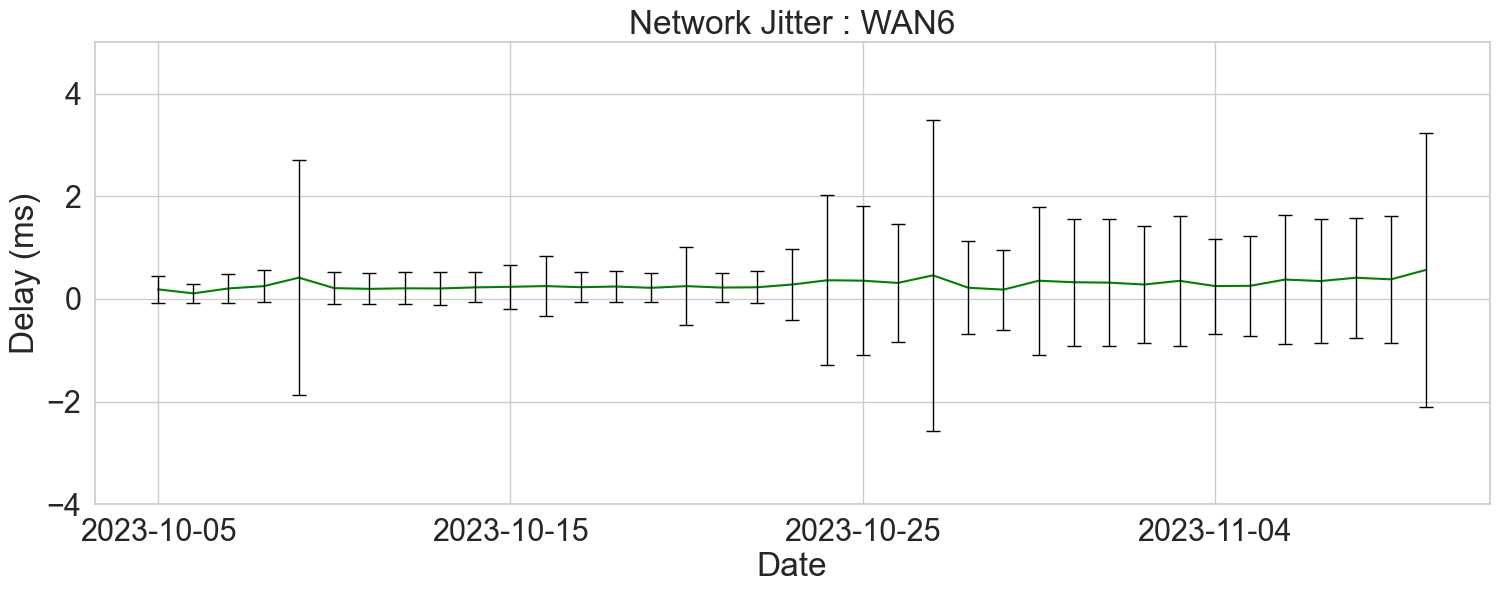} 
    \def\imageG{images/RIPE_LEO_JITTER.png} 
    \def\imageH{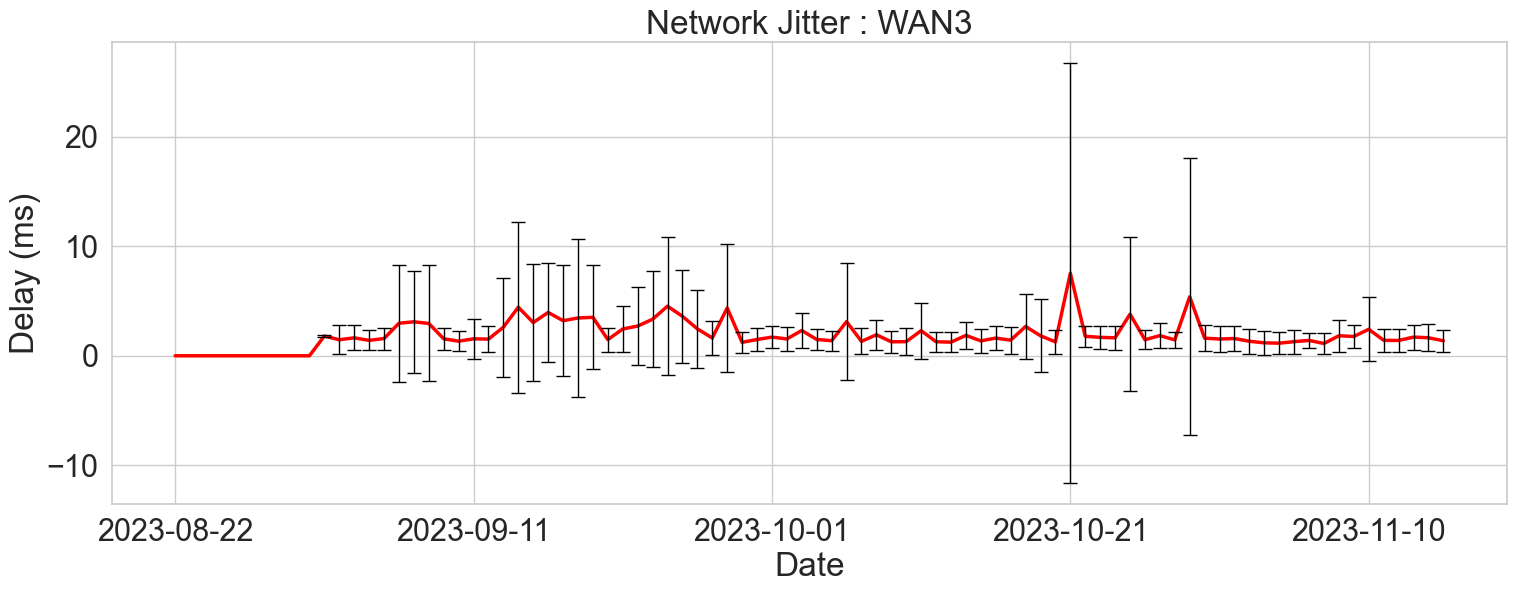} 
    \def\imageI{images/RIPE_DTN1_JITTER.png} 
    \def\imageJ{images/LOCAL_LEO_JITTER.png} 
    \def\imageK{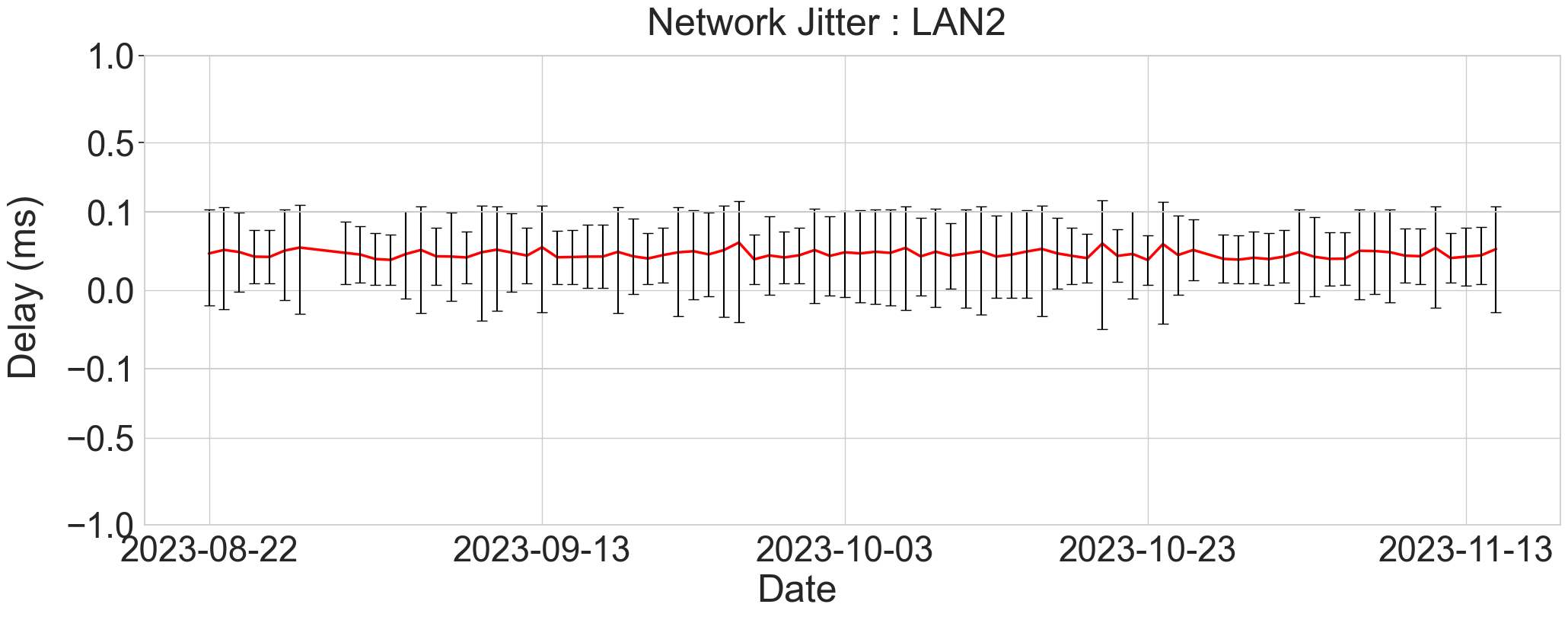} 
    \def\imageL{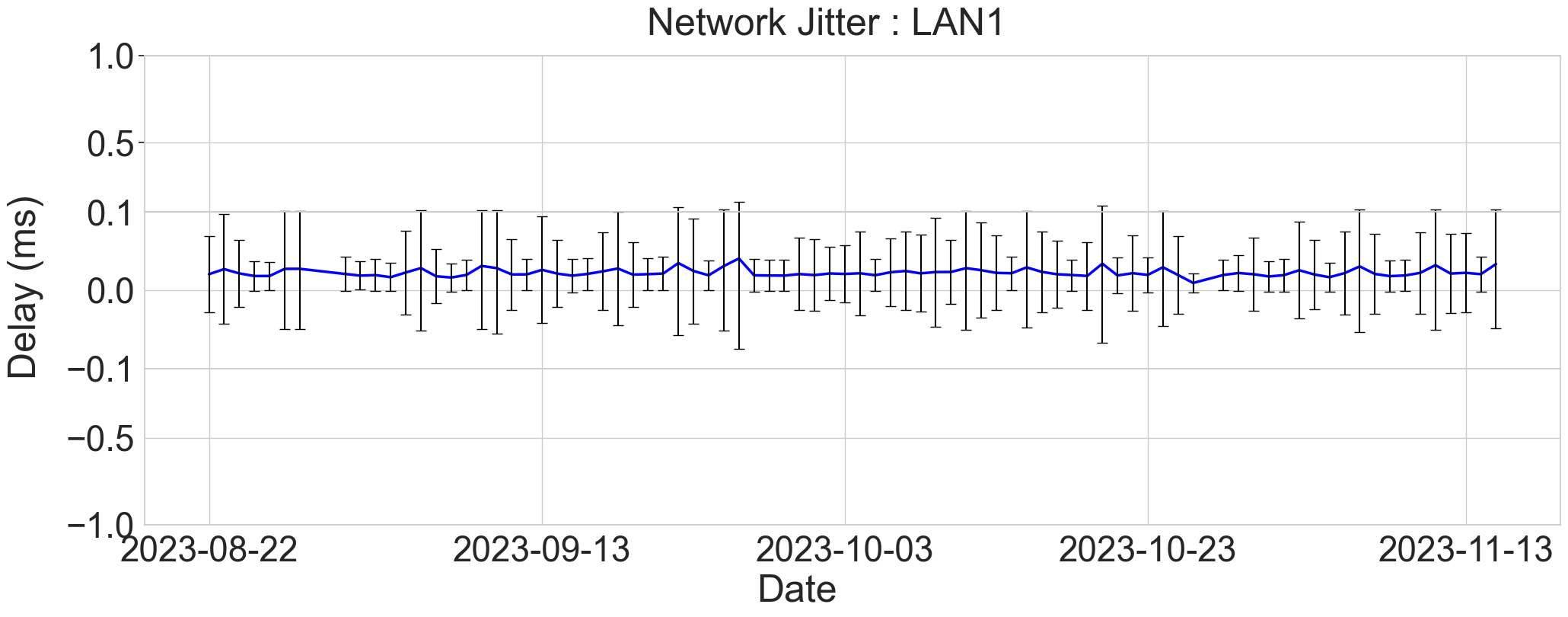} 
    
\begin{subfigure}[t]{0.3\textwidth}
        \includegraphics[width=\columnwidth]{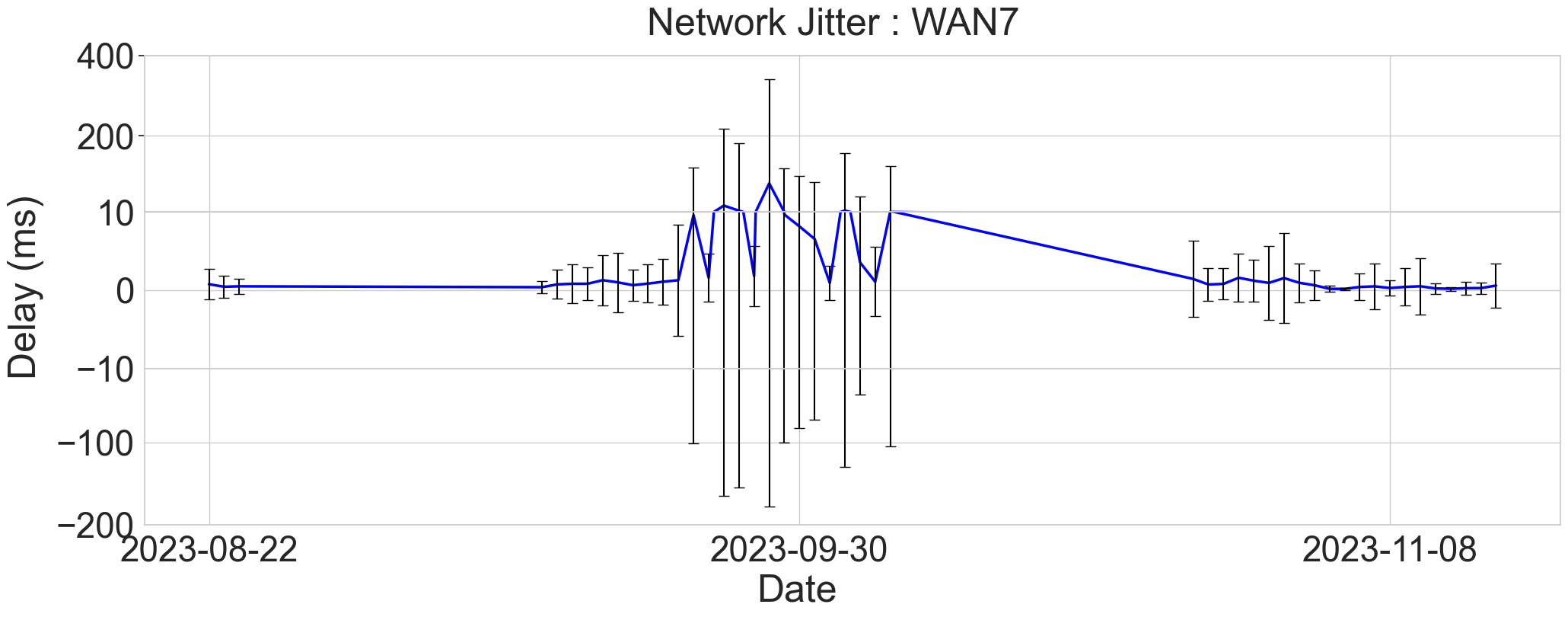}
        \caption{Jitter observed from Leo to external perfSONAR nodes}
        \label{fig:wan7_jitter}
\end{subfigure}
\begin{subfigure}[t]{0.3\textwidth}
        \includegraphics[width=\columnwidth]{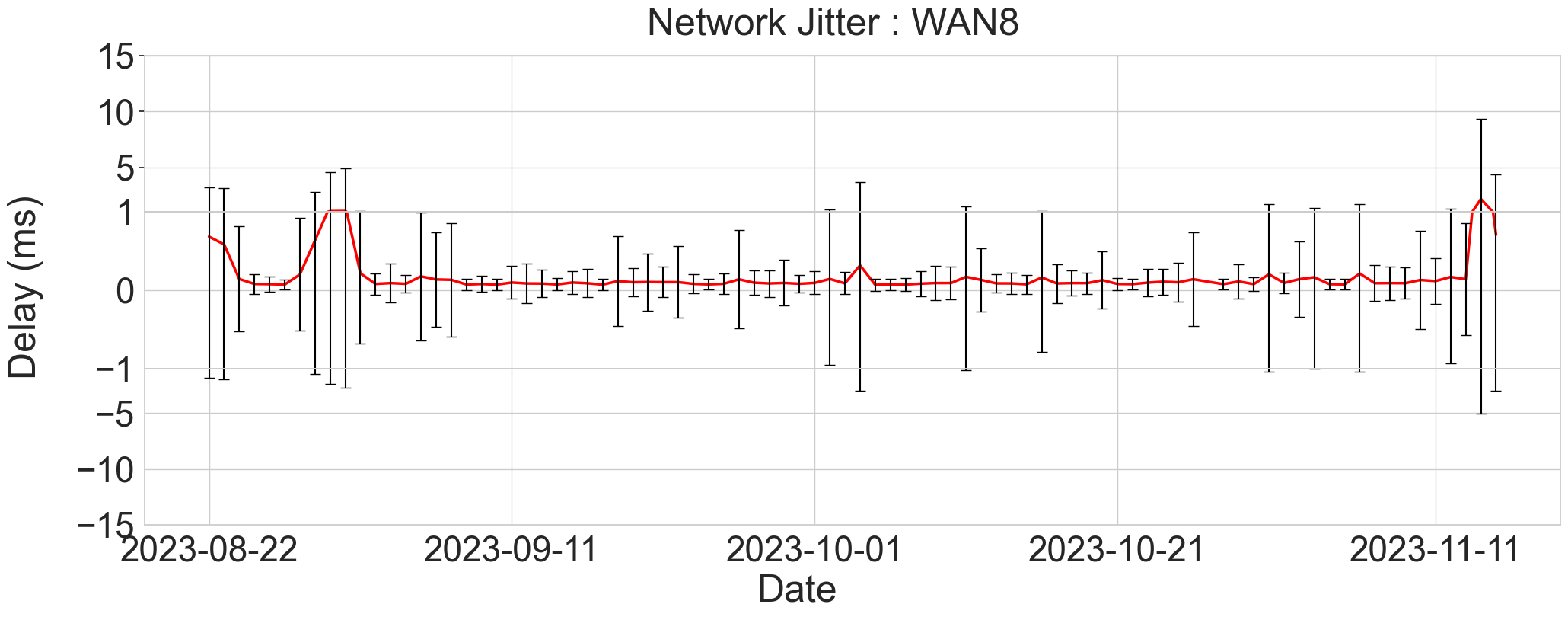}
        \caption{Jitter observed from perfSONAR1 to external perfSONAR nodes}
        \label{fig:wan8_jitter}
\end{subfigure}    
\begin{subfigure}[t]{0.3\textwidth}
        \includegraphics[width=\columnwidth]{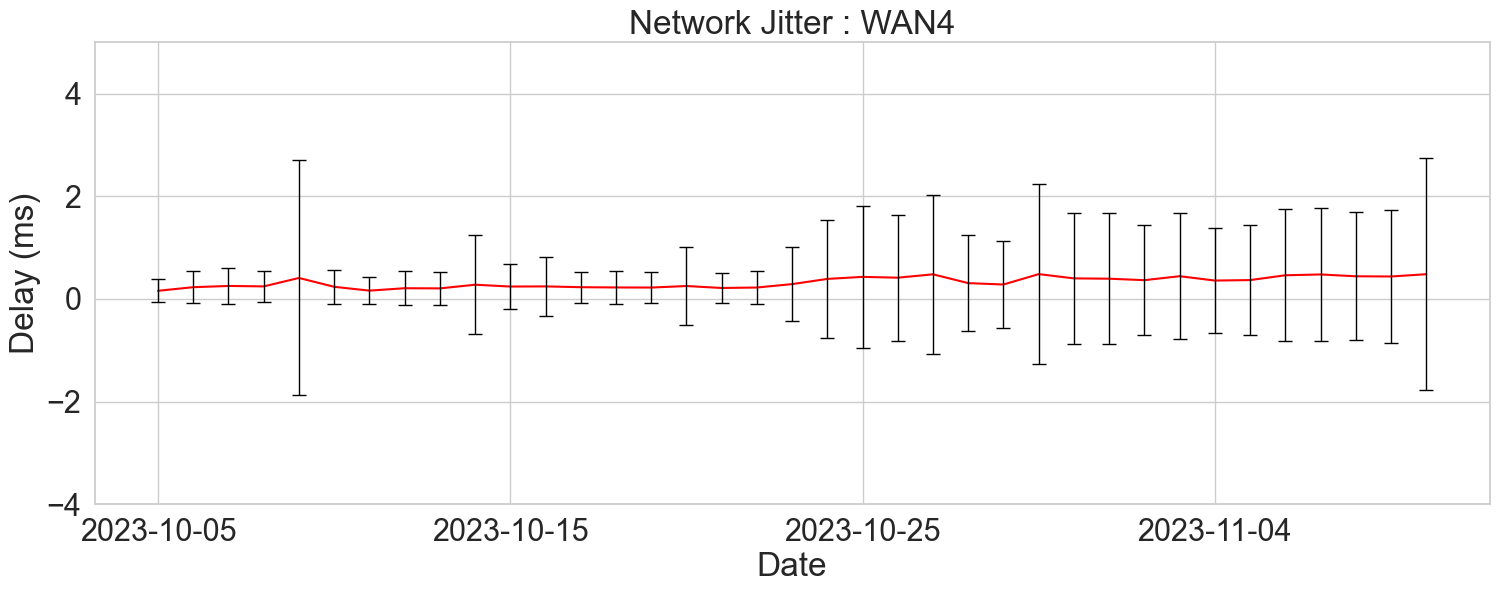}
        \caption{Jitter observed from GCP to perfSONAR1}
        \label{fig:wan4_jitter}
\end{subfigure}
\begin{subfigure}[t]{0.3\textwidth}
        \includegraphics[width=\columnwidth]{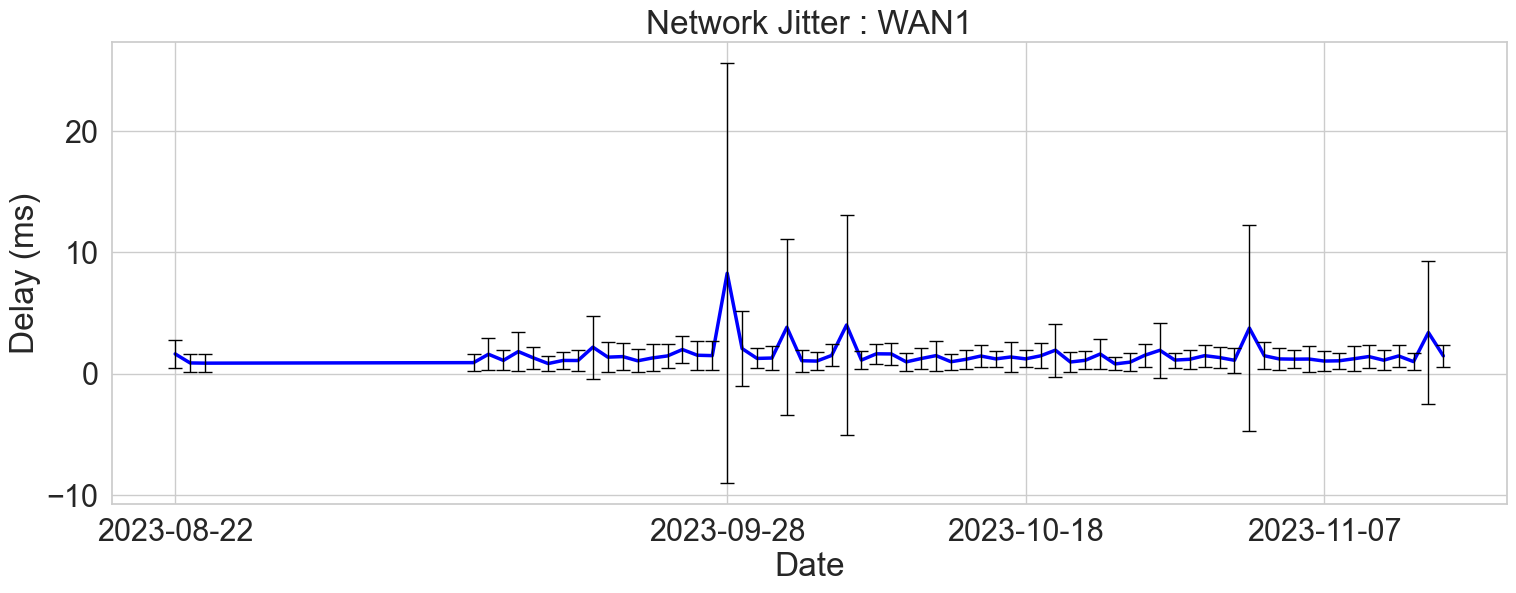}
        \caption{Jitter observed from RIPE Atlas to Leo}
        \label{fig:wan1_jitter}
\end{subfigure}
\begin{subfigure}[t]{0.3\textwidth}
        \includegraphics[width=\columnwidth]{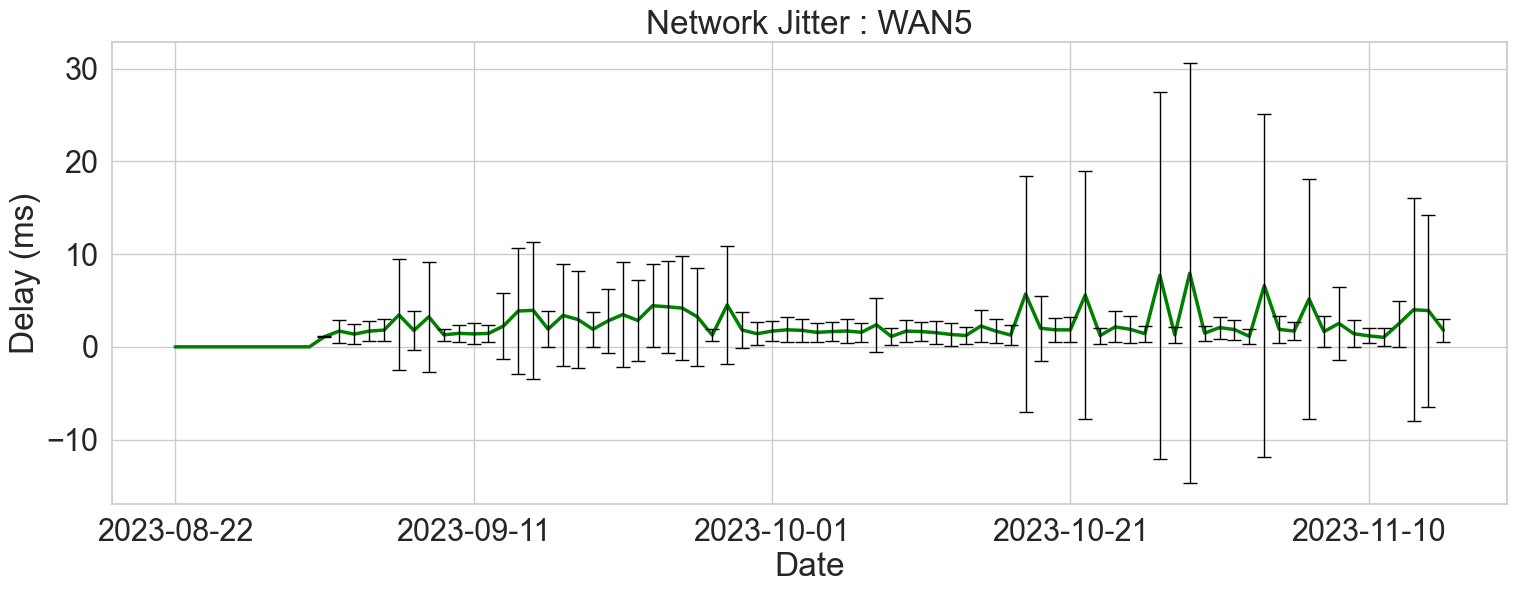}
        \caption{Jitter observed from one probe of RIPE Atlas to DTN1}
        \label{fig:wan5_jitter}
\end{subfigure}    
\begin{subfigure}[t]{0.3\textwidth}
        \includegraphics[width=\columnwidth]{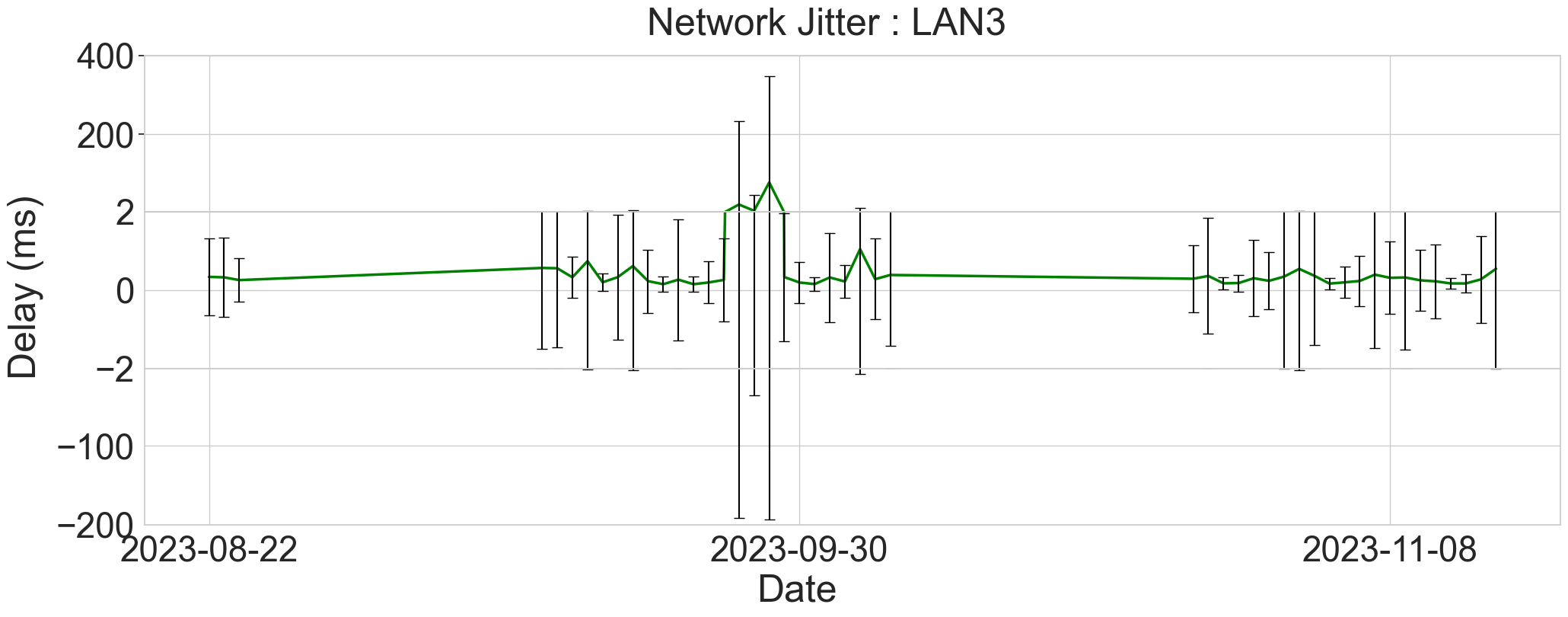}
        \caption{Jitter observed from Leo to perfSONAR1}
        \label{fig:lan3_jitter}
\end{subfigure}
\caption{Jitter Comparison}
    \label{fig:jitter}

\end{figure*}
Jitter is an important matrix for video and other real-time applications. In these experiments, we compare the jitter between the campus network and the DMZ.

When traffic is inbound from RIPE Atlas, Leo, the campus node, exhibits lower average jitter than the DTN1 or Perfsonar1 nodes as Figures \ref{fig:wan1_jitter} and \ref{fig:wan5_jitter} show. Jitter on the campus route tends to be 60-78\% lower than on the DMZ routes on average. DTN1 tends to exhibit higher variation in its daily jitter than Perfsonar1 by as much as 37 milliseconds, but the two DMZ nodes exhibit similar overall performance.

When traffic is inbound from GCP to Leo and the DMZ nodes, all three routes exhibit similar average jitter patterns between 0-1 milliseconds, only ever differing by fractions of milliseconds. Figure \ref{fig:wan4_jitter} represents the average jitter pattern, and differences in standard deviation from all three nodes are noted. Leo's route often experiences more deviation in its jitter than the DMZ nodes by as much as two milliseconds. The DTN1 node and Perfsonar1 node experienced a similar jitter pattern, so only the Perfsonar1 plot was selected to convey this experiment. However, the two nodes' difference in variation was noted. The DTN1 node experiences more deviation than the Perfsonar1 node by as much as 1.5 milliseconds.

    When traffic is outbound to perfSONAR nodes, as shown in Figures \ref{fig:wan7_jitter} and \ref{fig:wan8_jitter}, Leo and the DMZ nodes typically have an average jitter between 0-1 milliseconds. However, Leo often reaches higher jitter rates up to 10-63 milliseconds greater than the DMZ nodes. Both DTN1 and Perfsonar1 exhibited similar daily jitter and deviation patterns, so only the plot of Perfsonar1 was chosen to represent the DMZ for this experiment. 

    \textbf{Takeaways:}
    Campus networks experience more jitter than their DMZ counterparts. The average jitter on the campus network is also higher due to a higher number of competing flows. 
\subsection{Data Transfer Throughput}

One of the main reasons for creating DMZs is the higher data transfer rate that it enables. This section compares data transfer rates between the DMZ and the campus network. As mentioned earlier, for these tests, we downloaded publicly available Linux ISOs. We performed both experiments back to back to reduce variations in network conditions. Additionally, we did not tune the TCP stacks on the hosts. While such tuning significantly improves the data transfer rates, we wanted to establish a baseline comparison. Further tuning will improve data transfer performance in both DMZ and campus networks.

\begin{figure*}[!htbp]
    \centering
    \def\imageA{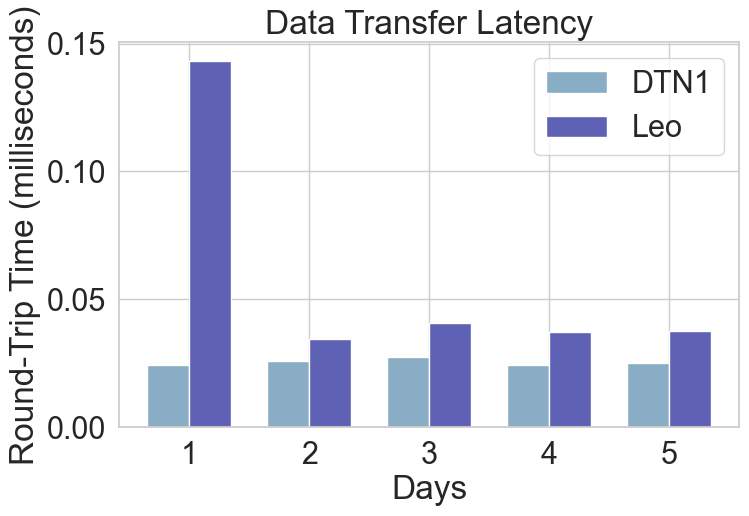}
    \def\imageB{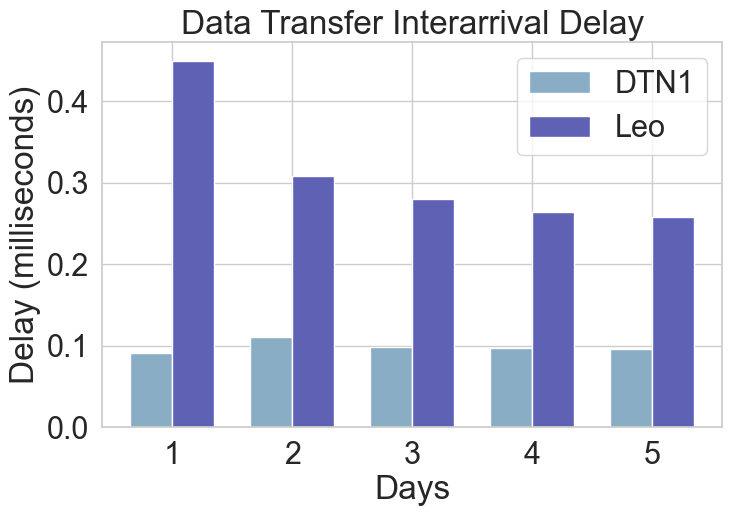}
    \def\imageC{images/DATA_TRANSFER_THROUGHPUT.png}
    \def\imageD{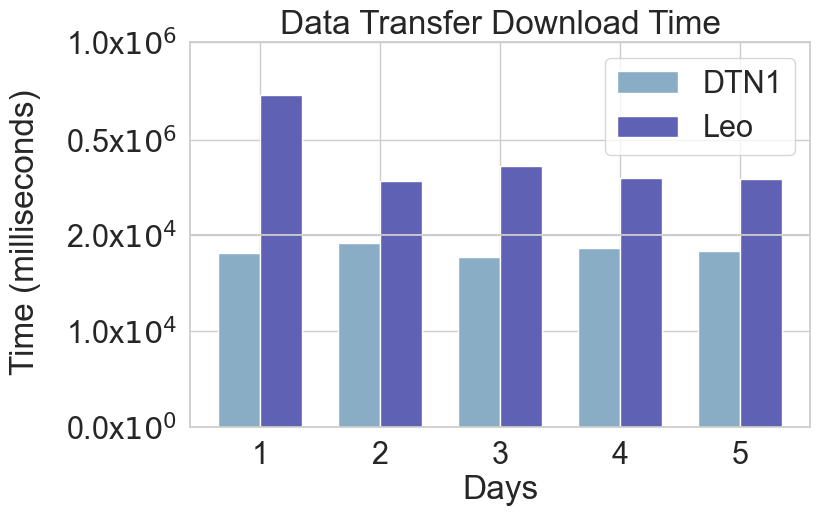} 
    \def\imageE{images/DATA_TRANSFER_BYTES_OUT.png} 
    \def\imageF{images/DATA_TRANSFER_RCVD_WIND.png} 
    \def\imageAA{images/UPLOAD_THROUGHPUT.png}
    \def\imageBA{images/BGP_HOP_COUNT.png}
    \def\imageBB{images/BGP_HOP_COUNT_2.png}
                
    \begin{subfigure}[t]{0.32\textwidth}
        \includegraphics[width=\columnwidth]{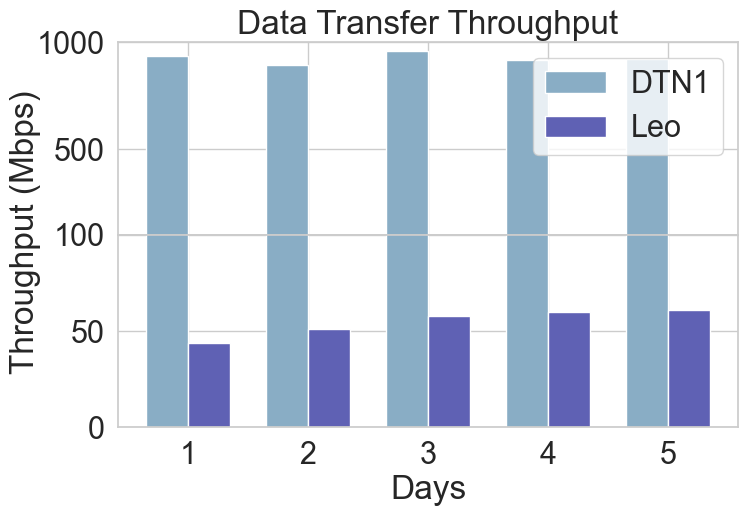}
        \caption{Average daily throughput observed from Data Transfer Experiments}
        \label{fig:datatf_throughput}
    \end{subfigure}
    \begin{subfigure}[t]{0.32\textwidth}
        \includegraphics[width=\columnwidth]{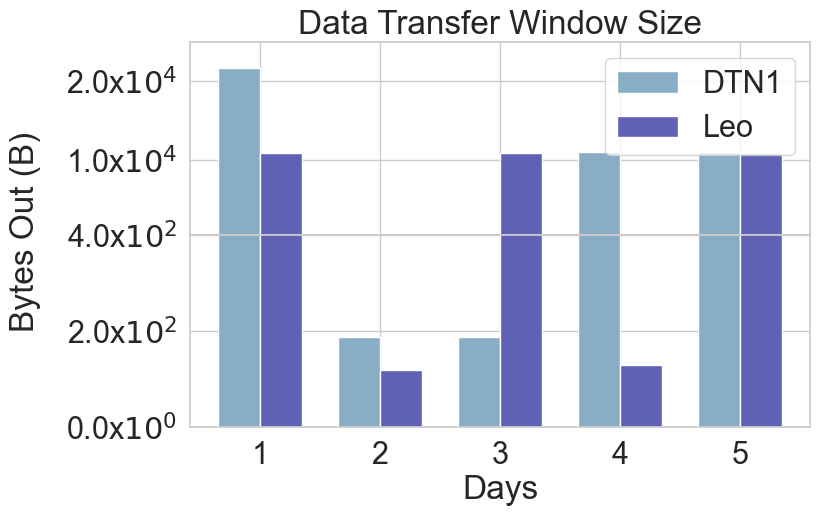}
        \caption{Average daily "bytes out" window size observed from Data Transfer Experiments}
        \label{fig:datatf_bytesout}
    \end{subfigure}
    \begin{subfigure}[t]{0.32\textwidth}
        \includegraphics[width=\columnwidth]{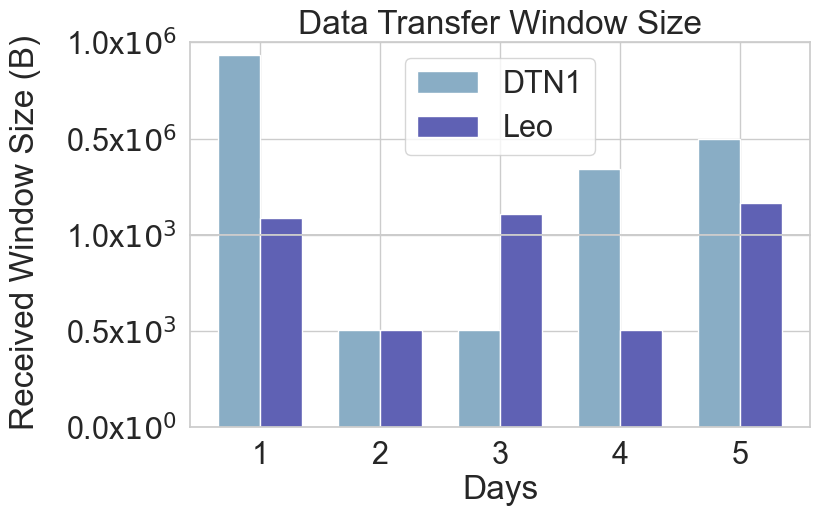}
        \caption{Average daily ``received" window size observed from Data Transfer Experiments}
        \label{fig:datatf_rcvd}
    \end{subfigure}
    \begin{subfigure}[t]{0.32\textwidth}
        \includegraphics[width=\linewidth]{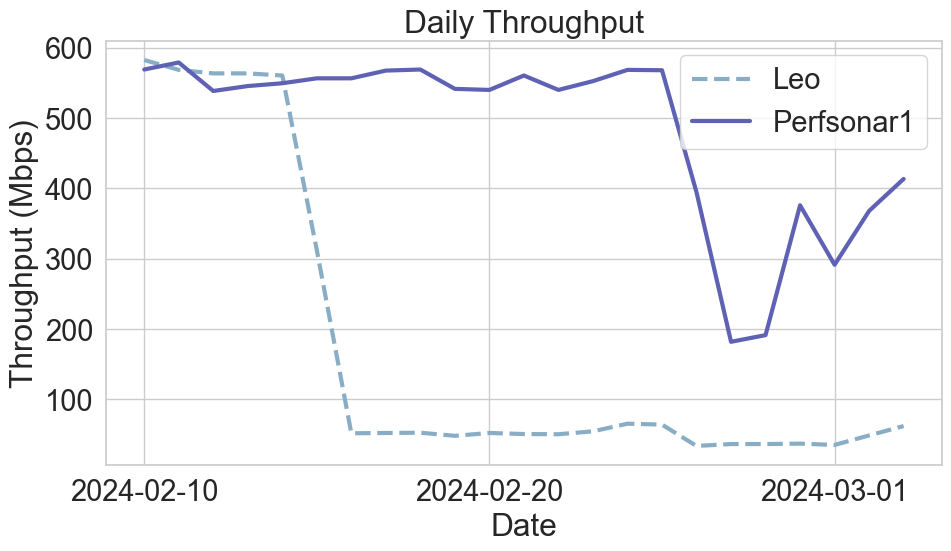}
        \caption{Average daily throughput observed from iperf3 experiments to GCP}
        \label{fig:iperf_upload}
    \end{subfigure}
    \begin{subfigure}[t]{0.32\textwidth}
        \includegraphics[width=\columnwidth]{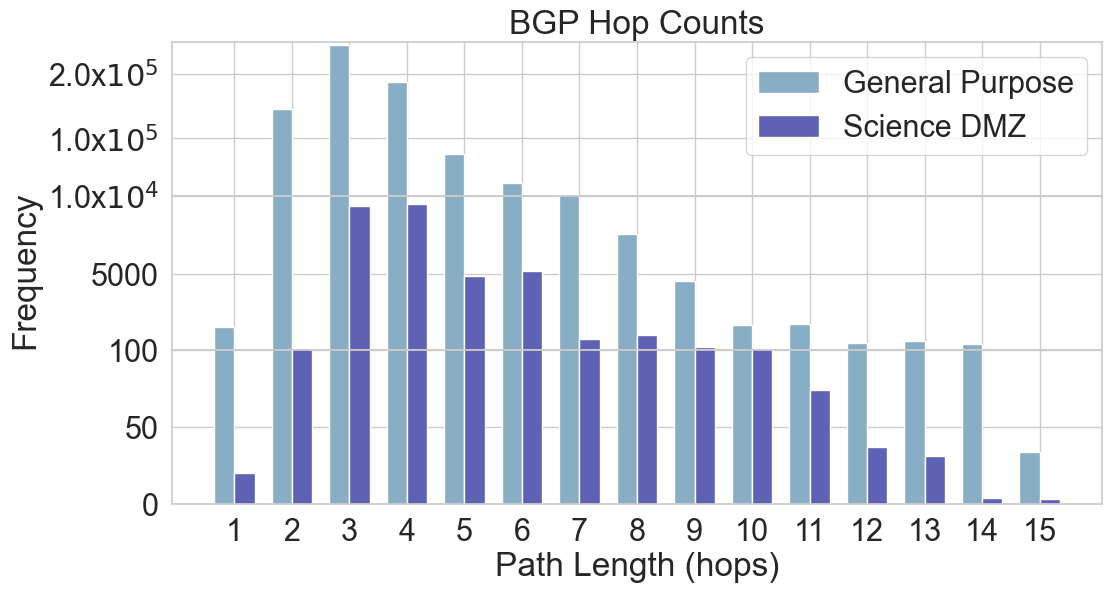}
        \caption{BGP Path Lengths}
        \label{fig:bgp1}
    \end{subfigure}
        \hfill
        \begin{subfigure}[t]{0.32\textwidth}
            \includegraphics[width=\columnwidth]{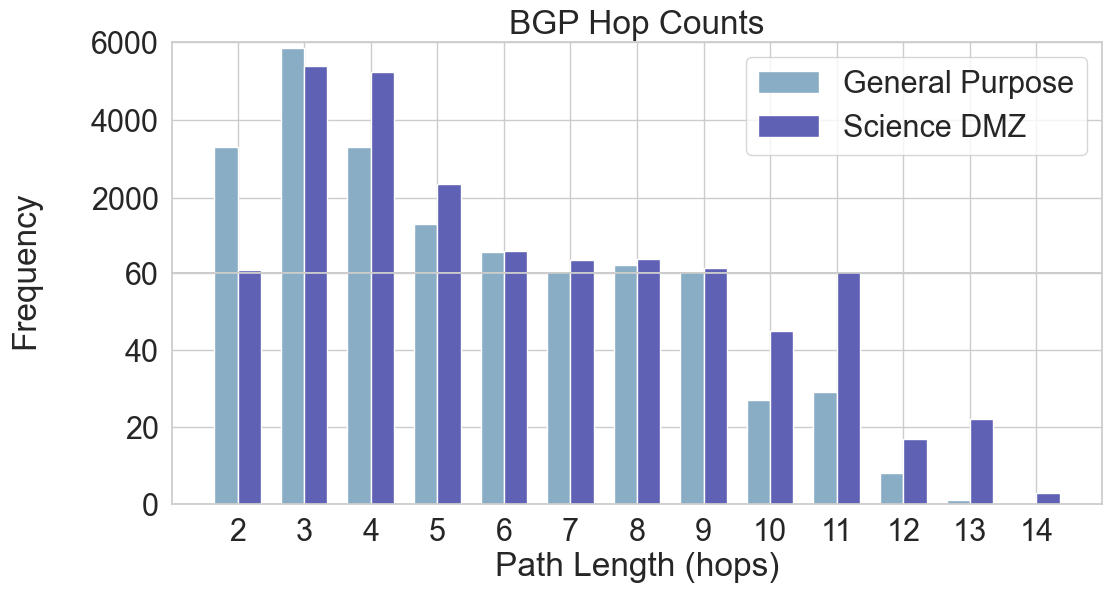}
            \caption{BGP Path Lengths}
            \label{fig:bgp2}
    \end{subfigure}
    \caption{Comparative Analysis of Data Transfer Metrics and BGP Path Lengths}
    \label{fig:data_transfer_analysis}
\end{figure*}

As Figure \ref{fig:datatf_throughput} shows, the average throughput was much higher on the DMZ when compared to the campus network. The host on the campus network could achieve only 50Mbps, while the host on the DMZ achieved close to 1Gbps. The slower data transfers are a result of packet loss and in-line firewall. On the other hand, the DMZ performs well since it only uses ACLs, and the loss rate is also low. 

We also looked at the TCP window sizes for these transfers, shown in Figures \ref{fig:datatf_bytesout} and \ref{fig:datatf_rcvd}. We looked at both the ``Bytes out" window size (bytes in flight) and the received window size, and the DTN had more oversized windows in both cases. The received window size was larger on Leo several times, but the throughput was low. This observation is consistent with what we would expect on a lossy link.

Figure \ref{fig:iperf_upload} corroborates these observations. We ran regular iperf3 tests between hosts on Google Cloud, DMZ, and the campus network. The DMZ host consistently outperforms the campus host in both upload and download performance.

\textbf{Takeaways}: The general purpose network performs significantly worse than a DMZ regarding file transfer performance. The campus network performs worse since it has more packet loss, the TCP window is smaller, and firewalls add latency to the packets.
\subsection{BGP Path Comparison}

This section compares the BGP path lengths between the DMZ and the campus network. We downloaded the BGP tables from the DMZ BGP router and campus ISP's BGP router.
First, we noticed that the commercial ISP had more additional routes than the Science DMZ router. The campus network had 715,810 BGP routes compared to 94,773 on the Science DMZ router. The campus BGP table also had three entries per destination as backup routes. We believe these are artifacts of BGP configurations. Other than having more route options in case of a failure and the capability of better load balancing, more BGP routes provide no additional advantages.

We then compared BGP hop counts between these networks. Figure \ref{fig:bgp1} shows the distribution. The general purpose network generally had a large number of paths with hop counts six or less (note the split Y axis). The DMZ also showed similar patterns.

Since the DMZ had less number of routes, we separated the intersection of these two tables and compared them in Figure \ref{fig:bgp2}. We found that the path lengths for the DMZ were slightly lower for shorter-length paths (hop counts \textless 3). For other DMZ routes, the hop count was larger than that of the campus routes. While BGP and IP path lengths are not always strictly correlated, these observations corroborate our findings in the previous experiments. 

\textbf{Takeaways:} The DMZ has less path diversity and longer path lengths than the campus network. While this may not directly affect performance, the resiliency of the DMZ can be improved by using additional fallback routes. Further, the path length can be reduced by creating better peering at the upstream, which requires negotiation with the upstream provider.




\section{Conclusions}
Science DMZs represent a paradigm shift in network design, tailored explicitly for scientific applications and distinct from traditional campus or general-purpose networks. The core principles of the Science DMZ, such as optimized paths for large data transfers and minimized security interference, position it as an advantageous environment for research and scientific collaboration. Over recent years, its adoption by numerous universities and organizations highlights its value in the academic and research communities.

Our comprehensive study over two years presents a nuanced picture. We confirm that the Science DMZ exhibits lower latency, higher throughput, and reduced jitter compared to general-purpose networks; when considering file transfer performance, the DMZ clearly outperforms general-purpose networks. The packet loss, smaller TCP windows, and added latency from firewalls in campus networks significantly hinder their efficiency in handling large-scale data transfers.

Science DMZs are not without their limitations. We observed non-intuitive results such as higher latency in specific scenarios and increased hop counts compared to campus networks. These findings suggest that while the Science DMZ can enhance certain aspects of network performance, it may not uniformly outperform campus networks in all areas, particularly in delay-sensitive applications like AR-VR. Our study reveals that the DMZ has less path diversity and longer path lengths than campus networks. While this impacts performance, strategic enhancements, such as developing better peering agreements and incorporating fallback routes, could mitigate these limitations.

In summary, while the Science DMZ model offers distinct advantages for specific research applications, it is not a one-size-fits-all solution. Such deployments must be carefully tailored to the particular needs and use cases of the communities they serve. This approach is paramount to fully harnessing the potential of DMZs in advancing scientific discovery and collaboration.

\bibliographystyle{ieeetr}
\bibliography{sample}

\begin{thebibliography}{10}

\bibitem{kissel2013efficient}
E.~Kissel, M.~Swany, B.~Tierney, and E.~Pouyoul, ``Efficient wide area data transfer protocols for 100 gbps networks and beyond,'' in {\em Proceedings of the Third International Workshop on Network-Aware Data Management}, NDM '13, (New York, NY, USA), Association for Computing Machinery, 2013.

\bibitem{dart2013science}
E.~Dart, L.~Rotman, B.~Tierney, M.~Hester, and J.~Zurawski, ``The science dmz: A network design pattern for data-intensive science,'' SC '13, (New York, NY, USA), Association for Computing Machinery, 2013.

\bibitem{lee2021scalable}
C.~Lee, M.~Jang, M.~Noh, and W.~Seok, ``Scalable design and algorithm for science dmz by considering the nature of research traffic,'' {\em The Journal of Supercomputing}, vol.~77, pp.~2979--2997, 2021.

\bibitem{nsfDMZ}
``{N}{S}{F} {A}ward {S}earch: {S}imple {S}earch {R}esults. --- nsf.gov.'' \url{https://www.nsf.gov/awardsearch/simpleSearchResult\\?queryText=dmz}.
\newblock [Accessed 13-03-2024].

\bibitem{8494729}
J.~Crichigno, E.~Bou-Harb, and N.~Ghani, ``A comprehensive tutorial on science dmz,'' {\em IEEE Communications Surveys \& Tutorials}, vol.~21, no.~2, pp.~2041--2078, 2019.

\bibitem{ESNET}
{The Energy Science Network}.

\bibitem{peisert2017medical}
S.~Peisert, E.~Dart, W.~Barnett, E.~Balas, J.~Cuff, R.~L. Grossman, A.~Berman, A.~Shankar, and B.~Tierney, ``{The medical science DMZ: a network design pattern for data-intensive medical science},'' {\em Journal of the American Medical Informatics Association}, vol.~25, pp.~267--274, 10 2017.

\bibitem{gonzalez2017bigdata}
A.~Gonzalez, J.~Leigh, S.~Peisert, B.~Tierney, E.~Balas, P.~Radulovic, and J.~M. Schopf, ``Big data and analysis of data transfers for international research networks using netsage,'' in {\em 2017 IEEE International Congress on Big Data (BigData Congress)}, pp.~344--351, 2017.

\bibitem{liu2017widearea}
Z.~Liu, P.~Balaprakash, R.~Kettimuthu, and I.~Foster, ``Explaining wide area data transfer performance,'' in {\em Proceedings of the 26th International Symposium on High-Performance Parallel and Distributed Computing}, HPDC '17, (New York, NY, USA), p.~167–178, Association for Computing Machinery, 2017.

\bibitem{caicedo2024machine}
C.~V. Caicedo, E.~F. Kfoury, J.~Gomez, J.~E. Pezoa, M.~Figueroa, and J.~Crichigno, ``Machine learning controller for data rate management in science dmz networks,'' {\em Computer Networks}, p.~110237, 2024.

\bibitem{vega2023explicit}
C.~F. Vega~Caicedo {\em et~al.}, ``Explicit feedback for congestion control in science dmz cyberinfrastructures based on programable data-plane switches.,'' 2023.

\bibitem{gegan2020anomaly}
R.~Gegan, C.~Mao, D.~Ghosal, M.~Bishop, and S.~Peisert, ``Anomaly detection for science dmzs using system performance data,'' in {\em 2020 International Conference on Computing, Networking and Communications (ICNC)}, pp.~492--496, IEEE, 2020.

\bibitem{mazloum2023enhancing}
A.~Mazloum, J.~Gomez, E.~Kfoury, and J.~Crichigno, ``Enhancing perfsonar measurement capabilities using p4 programmable data planes,'' in {\em Proceedings of the SC'23 Workshops of The International Conference on High Performance Computing, Network, Storage, and Analysis}, pp.~819--829, 2023.

\bibitem{crichigno2018comprehensive}
J.~Crichigno, E.~Bou-Harb, and N.~Ghani, ``A comprehensive tutorial on science dmz,'' {\em IEEE Communications Surveys \& Tutorials}, vol.~21, no.~2, pp.~2041--2078, 2018.

\bibitem{prasad6785344}
P.~Calyam, A.~Berryman, E.~Saule, H.~Subramoni, P.~Schopis, G.~Springer, U.~Catalyurek, and D.~K. Panda, ``Wide-area overlay networking to manage science dmz accelerated flows,'' in {\em 2014 International Conference on Computing, Networking and Communications (ICNC)}, pp.~269--275, 2014.

\bibitem{evalperf}
B.~Rababah, S.~Zhou, and M.~Bader, ``Evaluation the performance of dmz,'' {\em International Journal of Wireless and Microwave Technologies}, vol.~8, no.~1, pp.~1--13, 2018.

\bibitem{abhinit2022science}
I.~Abhinit, H.~Addleman, K.~Benninger, D.~DuRousseau, M.~Krenz, and B.~Meade, ``Science dmz: Secure high performance data transfer,'' tech. rep., 2022.

\bibitem{crichigno2021application}
J.~Crichigno, E.~Kfoury, E.~Bou-Harb, and N.~Ghani, ``Application and security aspects for large flows,'' in {\em High-Speed Networks: A Tutorial}, pp.~329--341, Springer, 2021.

\bibitem{crichigno2021data}
J.~Crichigno, E.~Kfoury, E.~Bou-Harb, and N.~Ghani, ``Data-link and network layer considerations for large data transfers,'' {\em High-Speed Networks: A Tutorial}, pp.~105--213, 2021.

\bibitem{RIPEatlas}
``Atlas console,'' 2022.

\bibitem{tierney2009perfsonar}
B.~Tierney, J.~Metzger, J.~Boote, E.~Boyd, A.~Brown, R.~Carlson, M.~Zekauskas, J.~Zurawski, M.~Swany, and M.~Grigoriev, ``perfsonar: Instantiating a global network measurement framework,'' {\em SOSP Wksp. Real Overlays and Distrib. Sys}, vol.~28, 2009.

\bibitem{ping}
``man page ping section 8.''

\bibitem{traceroute}
M.~Kerrisk, ``Traceroute(8) - linux manual page,'' 2021.

\bibitem{zurawski2013perfsonar}
J.~Zurawski, S.~Balasubramanian, A.~Brown, E.~Kissel, A.~Lake, M.~Swany, B.~Tierney, and M.~Zekauskas, ``perfsonar: On-board diagnostics for big data,'' in {\em Workshop on Big Data and Science: Infrastructure and Services}, Citeseer, 2013.

\bibitem{kratz2001ngi}
M.~Kratz, M.~Ackerman, T.~Hanss, and S.~Corbato, ``Ngi and internet2: Accelerating the creation of tomorrow's internet,'' in {\em MEDINFO 2001}, pp.~28--32, IOS Press, 2001.

\end{thebibliography}

\end{document}